\documentclass{aims} 
\usepackage{amsmath}
\usepackage{amssymb,amsfonts,amsthm,amsbsy,pifont,fleqn,euscript}
\usepackage{mathrsfs}
\usepackage{paralist}
\usepackage[misc]{ifsym}
\usepackage[colorlinks=true]{hyperref}
\hypersetup{urlcolor=blue, citecolor=red}
\allowdisplaybreaks

\def\currentvolume{X}
\def\currentissue{X}
\def\currentyear{200X}
\def\currentmonth{XX}
\def\ppages{X-XX}

\newtheorem{theorem}{Theorem}[section]

\newtheorem{lemma}[theorem]{Lemma}

\theoremstyle{definition}
\newtheorem{definition}[theorem]{Definition}
\newtheorem{remark}[theorem]{Remark}

\newcommand{\Om}{{\Omega}}

\newcommand{\Ac}{\mathfrak{A}}

\newcommand{\sE}{\mathscr{E}}

\newcommand{\R}{{\mathbb R}}

\newcommand{\cK}{{  K}}

\newcommand{\cW}{{  W}}

\newcommand{\wrt}{{with respect to }}

\numberwithin{equation}{section}
\allowdisplaybreaks[4]
\newcommand{\xqedhere}[2]{%
	\rlap{\hbox to#1{\hfil\llap{\ensuremath{#2}}}}}

\title[Global attractors for von Karman beam]
{Global attractors for a full von Karman beam transmission problem}
\author[Tamara Fastovska]{}



\begin{document}
	\maketitle

	\centerline{\scshape
		Tamara Fastovska $^{{\href{fastovskaya@karazin.ua}{\textrm{\Letter}}}1,2}$
	}
	
	\medskip
	
	{\footnotesize
		
		\centerline{$^1$Kharkiv Karazin National University, Ukraine}
	}
	
	\medskip
	
	{\footnotesize
		
		\centerline{$^2$Humboldt-Universit\"at zu Berlin, Germany}
	}
	
	\bigskip


	\begin{abstract}
		A  nonlinear transmisson problem for an elastic full von Karman beam
		is considered here.
		We prove that the system possesses a compact global attractor.
	\end{abstract}

	\section{Introduction} In this paper we consider a nonlinear
	transmission problem for an  elastic beam with the full von Karman nonlinearity.\@
	We assume that the beam, then in equilibrium,  occupies an
	interval $[0, L]$.\@ Let the part of the
	beam $(0, L_0)$, where $0<L_0<L$, is subjected to a structural damping
	while its complementary part $(L_0, L)$ is not.
	
	The system of differential equations for the transverse
	displacements $\phi(x,t)$, $u(x,t)$ and the longitudinal displacements $\omega(x,t)$,  $v(x,t)$ of the full von Karman part of the beam is as follows
	\begin{align}
		&\beta_1 \phi_{tt}-\mu_1 \phi_{ttxx}-\kappa\phi_{txx} +\lambda_1
		\phi_{xxxx}-\left(\left[\phi_x\left(\omega_x+1/2\phi_x^2\right) \right]_x \right)=g_1(x,t),\label{1}\\
		&\rho_1 \omega_{tt}+ \gamma \omega_{t}- \left(\omega_x+1/2\phi_x^2\right)_x=g_2(x,t).\qquad\qquad\qquad\;\;\;t>0,\;\;x\in (0, L_0)\label{2}\\
		&\beta_2 u_{tt}-\mu_2 u_{ttxx}+\lambda_2
		u_{xxxx}-\left(\left[u_x\left(v_x+1/2u_x^2\right) \right]_x \right)=g_3(x,t),\qquad\qquad\qquad\;\;\label{3}\\
		&\rho_2 v_{tt}- \left(v_x+1/2u_x^2\right)_x=g_4(x,t),\qquad\qquad\qquad\qquad\;\;\;\;\;\;t>0,\;\;x\in ( L_0, L)\label{4}
	\end{align}
	Here $\rho_i, \kappa, \beta_i, \lambda_i, \mu_i, \gamma$ for $i=1,2$ are positive constants.\par
	System \eqref{1}-\eqref{4} is supplemented with the transmission boundary conditions
	\begin{align}
		\phi(L_0, t)=u(L_0, t),& \;\omega(L_0, t)=v(L_0, t),\label{6}\\ \phi_x(L_0, t)=u_x(L_0, t),\;
		\lambda_1 \phi_{xx}(L_0, t)&= \lambda_2 u_{xx}(L_0, t),
		\;\omega_x(L_0, t) = v_x(L_0, t),\label{7}\\
		(\lambda_1 \phi_{xxx}-\mu_1 \phi_{ttx}-\kappa \phi_{tx}) (L_0,
		t)&=(\lambda_2 u_{xxx}-\mu_2 u_{ttx}) (L_0,
		t),\label{8}
	\end{align}
	and boundary conditions on the ends of the beam
	\begin{align}
		\phi_{x}(0, t)=0,\; \phi(0, t)=0,\;\omega(0, t)=0,\; u(L, t)=0,\; u_{xx}(L, t)=0,\; v(L, t)=0.\label{18}
	\end{align}
	We impose also initial conditions
	\begin{align}
		&\phi(x,0)=\phi_0(x), \;\phi_t(x,0)=\phi_1(x),\;\omega(x,0)=\omega_0(x), \;\omega_t(x,0)=\omega_1(x),\\
		&u(x,0)=u_0(x),  \;u_t(x,0)=u_1(x),\;v(x,0)=v_0(x),  \;v_t(x,0)=v_1(x)
		, \;x\in \Omega.\label{24}
	\end{align}
	Problems related to the study of stabilization and long-time
	behavior of transmission problems have  attracted an ample attention.  Several works are
	devoted to the transmission problems for Kirchhoff thermoelasticity
	(see, e.g. \cite{PO1, Vila, Potomkin}). Paper \cite{PO1} is related to the linear transmission problem between elastic and thermoelastic  Kirchhoff
	beams with the classical Fourier law of heat conduction. In
	\cite{Vila} the linear Kirchhoff problem with localized thermal dissipation
	of hereditary type is considered. In both works the exponential decay rate of the energy is shown. In paper    \cite{Potomkin} a nonlinear transmission problem for elastic and thermoelastic plates is under consideration. In case of  the Berger type nonlinearities, the existence of a compact global attractor is established.  Paper \cite{F1} is devoted to a linear transmission problem for a Kirchhoff-Timoshenko beam with different types of heat conduction, the exponential stability of the system is established. \par
	In the present paper we investigate the long-time
	dynamics of an elastic  beam  described by the full von Karman model whose part is subjected to a structural damping.
	\par
	
	The main aim of the present paper is to investigate the
	asymptotic behavior of the solutions to the problem considered no
	matter how small the dissipative  part of the beam
	is.  The complex structure of the nonlinear terms does not allow to prove the existence of an absorbing ball directly. To overcome this difficulty and to show the existence of a compact global attractor we establish the gradient property by means of a unique continuation result. \par
	The paper is organized as follows.  In
	Section 2 we formulate standard results on the existence of global attractors, introduce notations and state a  well-posedness theorem for the problem considered.
	Section 3 is devoted to the asymptotic compactness  of the system. In Section 4 the main result on the existence of a global attractor is established.\par
	
	\section{Preliminaries, notations, and well-posedness.}
	\subsection{Abstract results on  attractors.}
	For the readers' convenience  we recall some basic definitions and results from the theory of attractors.
	\begin{definition}[\cite{BabinVishik,chDQ,Chueshov,CL,Temam}]
		A global attractor of a  dynamical system  $(S_t, H)$ with the evolution operator $S_t$ on a complete metric space $H$ is defined as a bounded closed  set $\Ac\subset H$ which is  invariant ($S_t\Ac=\Ac$ for all $t>0$) and  uniformly attracts all other bounded  sets:
		\begin{equation*}
			\lim_{t\to\infty} \sup\{{\rm dist}_H  (S_ty,\Ac):\ y\in B\} = 0
			\quad\mbox{for any bounded  set $B$ in $H$.}
		\end{equation*}
	\end{definition}
	To establish the existence of attractor we use the concept of gradient systems. The main feature of these systems
	is that in the proof of the existence of a global attractor we can avoid a dissipativity
	property (existence of an absorbing ball) in the explicit form (\cite{chDQ}).
	\begin{definition}[\cite{chDQ,Chueshov,CFR,CL}]
		Let $Y \subseteq H$ be a forward invariant set of a dynamical system
		$(S_t, H)$. A continuous functional $L(y)$ defined on $Y$ is said to be a Lyapunov function on
		$Y$ for the dynamical system $(S_t, H)$ if $t\mapsto L(S_ty)$ is a nonincreasing function for
		any $y\in Y$.\par
		The Lyapunov function  is said to be strict on $Y$ if the equation $L(S_ty) =
		L(y)$ for all $t > 0$ and for some $y \in Y$ implies that $S_ty = y$ for all $t > 0$; that is, $y$
		is a stationary point of $(S_t, H)$.\par
		The dynamical system  is said to be gradient if there exists a strict Lyapunov
		function on the whole phase space $H$.
	\end{definition}
	\begin{definition}[\cite{chDQ,Chueshov,CFR,CL}]
		A dynamical system $(X,S_t)$ is said to be asymptotically smooth
		if for any  closed bounded set $B\subset X$ that is positively invariant ($S_tB\subseteq B$)
		one can find a compact set $\cK=\cK(B)$ which uniformly attracts $B$, i.~e.
		$\sup\{{\rm dist}_X(S_ty,\cK):\ y\in B\}\to 0$ as $t\to\infty$.
	\end{definition}
	In order to prove  the asymptotical smoothness of  system
	(\ref{1})-(\ref{24}) we rely on the  compactness criterion due to
	\cite{Khanmamedov}, which is recalled below in an abstract version
	formulated in \cite{CL}.
	\begin{theorem}{\cite{CL}} \label{theoremCL} Let $(S_t, H)$ be a dynamical system on a complete metric
		space $H$ endowed with a metric $d$. Assume that for any bounded positively invariant
		set $B$ in $H$ and for any $\varepsilon>0$ there exists $T = T (\varepsilon, B)$ such that
		\begin{equation}
			\label{te}
			d(S_T y_1, S_T y_2) \le \varepsilon+ \Psi_{\varepsilon,B,T} (y_1, y_2), y_i \in B ,
		\end{equation}
		where $\Psi_{\varepsilon,B,T} (y_1, y_2)$ is a function defined on $B \times B$ such that
		\[
		\liminf\limits_{m\to\infty}\liminf\limits_{n\to\infty}\Psi_{\varepsilon,B,T} (y_1, y_2) = 0
		\]
		for every sequence ${y_n} \in B$. Then $(S_t, H)$ is an asymptotically smooth dynamical
		system.
	\end{theorem}
	The following statement collects criteria on existence and properties of attractors to gradient systems.
	\begin{theorem}[\cite{chDQ, CFR, CL}]
		\label{abs}
		Assume that $(S_t, H)$ is a gradient asymptotically smooth  dynamical
		system. Assume its Lyapunov function $L(y)$ is bounded from above on any
		bounded subset of $H$ and the set $\cW_R=\{y: L(y) \le R\}$ is bounded for every $R$. If the
		set $\EuScript N$ of stationary points of $(S_t, H)$ is bounded, then $(S_t, H)$ possesses a compact
		global attractor. Moreover,
		the global attractor $\Ac$ consists of  full trajectories
		$\gamma=\{ U(t)\, :\, t\in\R\}$ such that
		\begin{equation}\label{conv-N}
			\lim_{t\to -\infty}{\rm dist}_{H}(U(t),\EuScript N)=0 ~~
			\mbox{and} ~~ \lim_{t\to +\infty}{\rm dist}_{H}(U(t),\EuScript N)=0.
		\end{equation}
		and
		\begin{equation}\label{7.4.1}
			\lim_{t\to +\infty}{\rm dist}_{H}(S_tx,\EuScript N)=0
			~~\mbox{for any $x \in H$;}
		\end{equation}
		that is, any trajectory stabilizes to the set $\EuScript N$ of  stationary points.
	\end{theorem}
	\subsection{Notations.}
	Let $D$ be a bounded interval in $\mathbb R$  and $s\in\mathbb R$.
	We denote  by $H^s(D)$ the standard Sobolev space of order $s$
	on a set $D$ which we define as restriction (in the sense of distributions)
	of the
	space $H^s(\mathbb R)$ (introduced via Fourier transform).
	We denote by $\|\cdot \|_{s}$ the norm in  $H^s(D)$
	which we define by the relation
	$
	\|f\|_{s}^2=\inf\left\{\|g\|_{s,\mathbb R}^2\, :\; g\in H^s(\mathbb R),~~ g=f ~~
	\rm{on}~~D
	\right\}.$
	We also use the notation $\|\cdot \|=\|\cdot \|_{0}$
	for the corresponding $L_2$-norm and, similarly, $(\cdot,\cdot)$ for the $L_2$
	inner product.
	We denote by $H^s_0(D)$ the closure of $C_0^\infty(D)$ in $H^s(D)$
	(\wrt  $\|\cdot \|_{s}$) and by $H^s_{\{M\}}(D)$ the closure of $\{f\in C^\infty(D): f(M)=0\}$, where $M$ can denote $0$, $L$, or $L_0$.\par
	For the  component $\xi=(\phi, u)$ we  define the space
	\begin{equation*}\label{X-space}
		X=
		\{(\phi, u)\in H^2_{\{0\}}(0, L_0) \times (H^2\cap H^1_{\{L\}})(L_0, L): \phi(L_0)=u(L_0),  \phi_x(L_0)=u_x(L_0)\}.
	\end{equation*}
	We also define a  space for the  component $\zeta=(\omega,v)$
	\begin{equation*}\label{Y-space}
		Y=\{(\omega,v)\in H^1_{\{0\}}(0, L_0) \times H^1_{\{L\}}(L_0, L):\omega(L_0)=v(L_0) \}.
	\end{equation*}
	We equip the space for $Z=(\xi, \zeta)$
	\begin{equation}
		\label{V-space}
		V=X\times Y
	\end{equation}
	with the inner  product
	$$(Z_1, Z_2)_V\!=\!\int_0^{L_0} \omega_{1x}\omega_
	{2x} dx \!+\!\lambda_1 \int_0^{L_0} \phi_{1xx} \phi_{2xx} dx\!+\!\int_{L_0}^L v_{1x}v_
	{2x} dx \!+\!\lambda_2 \int_{L_0}^L u_{1xx} u_{2xx} dx,
	$$
	where $Z_i=(\xi_i, \zeta_i)$, $i=1,2$.\par
	We also define the space
	\begin{equation*}\label{Y-space1}
		\tilde Y=L_2(0, L_0)\times L_2(L_0, L),
	\end{equation*}
	and
	\begin{equation}
		\label{W-space}
		W= Y\times \tilde Y
	\end{equation}
	endowed with the inner product
	\begin{align*}
		(Z_1, Z_2)_{W}&=\beta_1 \int_0^{L_0}  \phi_1 \phi_2 dx+\mu_1  \int_0^{L_0} \phi_{1x}  \phi_{2x} dx+\rho_1 \int_{L_0}^L \omega_1  \omega_2 dx\\ &\quad+\beta_2 \int_{L_0}^L  u_1 u_2 dx+\mu_2  \int_{L_0}^L  u_{1x}  u_{2x} dx+\rho_2 \int_{L_0}^L v_1  v_2 dx.
	\end{align*}
	As the phase space we use
	\begin{equation}
		\label{h} H=V \times W.
	\end{equation}
	Throughout the paper we will denote by $C$ a generic positive constant.
	\subsection{Well-posedness.}
	To show the well-posedness of problem \eqref{1}-\eqref{24} we will need the following auxiliary result, which follows straightforward from embedding theorems (see, e.g. \cite{BabinVishik}).
	\begin{lemma}
		\label{lem1}
		There exists a positive constant $C$ such that  for any
		$z=(u; v)\in  H^2(0,L)\times H^1(0,L)$
		we have that
		\begin{equation*}
			\|z\|_{H^1(0,L)}^2\le C\left(Q(\zeta)+\|u\|_{H^{2}(0,L)}^4\right).
		\end{equation*}
		Here $Q(\xi)=\int_\Omega (v_x+\frac{u_x^2}{2})^2 dx$.
	\end{lemma}
	We will also use  the following lemma.
	\begin{lemma}
		\label{lem2}
		There exists a positive constant $C$ such that for any
		$g\in H^{1}(0,L)$ we have that
		\begin{equation*}
			\max\limits_\Om|g|\le C\|g\|_{H^{1/2+\delta}(0,L)},\; \text{for any}\; \delta>0.
		\end{equation*}
		
	\end{lemma}
	The ideas of the proof can be found e.g. in \cite{Boutet}.\par
	We define the spaces of test functions
	$$
	\EuScript L_T=\left\{\Psi=(\Psi_1, \Psi_2, \Psi_3, \Psi_4):
	\Psi\in L_2(0,T; V), \Psi_t\in
	L_2(0,T; W)\right\}
	$$
	and
	$\EuScript L_T^0=\{\Psi\in \EuScript L_T\, :\, \Psi(T)=0\}$.
	\par
	We also define positive self-adjoint operators
	\begin{equation*}
		N_1(\phi,  u)=(\lambda_1 \phi_{xxxx},  \lambda_2 u_{xxxx}): \EuScript D(N_1)\subset \tilde Y\mapsto \tilde Y
	\end{equation*}
	and
	\begin{equation*}
		N_2(\omega, v)=( -\omega_{xx},  -v_{xx}): \EuScript D(N_2)\subset \tilde Y\mapsto \tilde Y
	\end{equation*}
	with the domains
	\begin{multline*}
		\EuScript D(N_1)=\{(\phi, u)\in X\cap (H^4
		(0,L_0)\times H^4(L_0,L)):\lambda_1 \phi_{xx} (L_0,
		t)=\lambda_2 u_{xx} (L_0,
		t),\\ \lambda_1 \phi_{xxx} (L_0,
		t)=\lambda_2 u_{xxx} (L_0,
		t),\;u_{xx} (L,
		t)=0\}
	\end{multline*}
	and
	\begin{equation*}
		\label{bl}
		\EuScript D(N_2)=\{(\omega, v)\in Y\cap (H^2
		(0,L_0)\times  H^2
		(L_0,L)):\omega_x(L_0, t)= v_x(L_0, t)\}.
	\end{equation*}
	We also introduce a bounded operator $G: Y\mapsto Y'$ as follows
	\begin{align}&(G(\phi_1, u_1), (\phi_2, u_2))_{L_2(0, L_0)\times L_2(L_0, L)} \nonumber\\&= \beta_1 \int_0^{L_0}  \phi_1 \phi_2 dx+\mu_1  \int_0^{L_0} \phi_{1x}  \phi_{2x} dx+\beta_2 \int_{L_0}^L  u_1 u_2 dx+\mu_2  \int_{L_0}^L  u_{1x}  u_{2x} dx.
	\end{align}
	It is easy to see that $G$ is an isomorphism of $Y$ onto $Y'$. Let us consider the operator $G^{-1} N_1:\EuScript D(G^{-1}N_1)\subset Y\mapsto Y$.\par
	In order to make our statements precise we need to introduce the definition of weak solutions to problem (\ref{1})-(\ref{24}).
	\begin{definition}
		A  function $Z(t)=(\xi(t), \zeta(t))$, where $\xi(t)=(\phi(t), u(t))$ and $\zeta(t)=(\omega(t), v(t))$ is said to be a weak solution to
		problem (\ref{1})-(\ref{24}) on a time interval $[0,T]$ if
		\begin{itemize}
			\item $Z\in L_\infty(0,T;V)$, $Z_t\in L_\infty(0,T;W)$;
			\item $Z(0)=Z_0=(\phi_0, u_0,\omega_0, v_0)$ ;
			\item for every $\Psi=(\Psi_1, \Psi_2, \Psi_3, \Psi_4)\in
			\EuScript L_T^0$  the following equality holds
			\begin{align}
				\label{sol_def}
				&-\int_0^T \int_0^{L_0}(\beta_1\phi_{t} \Psi_{1t}+\rho_1\omega_{t} \Psi_{3t})dx dt-
				\int_0^T \int_{L_0}^L(\beta_2 u_{t} \Psi_{2t}+\rho_2 v_{t} \Psi_{4t})dx dt \nonumber\\& \quad-\mu_1\int_0^T \int_0^{L_0} \phi_{tx} \Psi_{1tx} dx dt-\mu_2\int_0^T \int_{L_0}^L u_{tx} \Psi_{2tx} dx dt+\kappa\int_0^T \int_0^{L_0} \phi_{tx} \Psi_{1x} dx dt \nonumber\\&\quad+\int_0^T\int_0^{L_0}  K_1(Z, \Psi)dx dt+ \int_0^T\int_{L_0}^{L}  K_2(Z, \Psi)dx dt+\lambda_1 \int_0^T\int_0^{L_0}  \phi_{xx} \Psi_{1xx}dx dt\nonumber\\& \quad+\lambda_2 \int_0^T\int_{L_0}^{L}  u_{xx} \Psi_{2xx}dx dt
				= \int_0^{L_0}(\beta_1\phi_{1} \Psi_{1}(0)+\rho_1\omega_{1} \Psi_{3}(0))dx \nonumber\\ &\quad+ \int_{L_0}^L(\beta_2 u_{1} \Psi_{2}(0)+\rho_2v_{1} \Psi_{4}(0))dx+\mu_1 \int_0^{L_0} \phi_{1x} \Psi_{1x}(0) dx+\mu_2 \int_{L_0}^L u_{1x} \Psi_{2x}(0) dx \nonumber\\ &\quad+\int_0^T( \sum_{i=1, 2}\int_0^{L_0}g_i\Psi_i  dx +\sum_{i=3,4}\int_{L_0}^L g_i\Psi_i  dx) dt,
			\end{align}
			where
			\begin{equation*}\label{def_r}
				K_1(Z, \Psi)\!=\!\int_0^{L_0}\!\!(\omega_x+\frac{\phi_x^2}{2}) (\Psi_{1x}+\phi_x\Psi_{3x})dx,\,K_2(Z, \Psi)\!=\!\int_{L_0}^L\!(v_x+\frac{u_x^2}{2}) (\Psi_{2x}+u_x\Psi_{4x})dx.
			\end{equation*}
		\end{itemize}
	\end{definition}
	
	The well-posedness result is as follows.
	\begin{theorem} \label{th:WP}
		Assume that
		\begin{equation}\label{g}g(x, t)=(g_1, g_2, g_3, g_4)\in L_2(0,T; W'),\end{equation} $$ U_0=(Z_0, Z_1)=(\phi_0, u_0, \omega_0, v_0, \phi_1, u_1, \omega_1, v_1)\in H.$$   Then
		for any interval $[0,T]$
		there exists a unique weak solution $Z(t)$ to (\ref{1})-(\ref{24})
		with the initial data $U_0$. This solution possesses the following properties:
		\begin{enumerate}
			\item[(i)]  $U(t, U_0)=(Z(t); Z_t(t))\in C(0,T; H).$
			\item[(ii)] The solution depends continuously  on initial data, i.e.
			if $U_n\to U_0$ in the norm of $H$, then  $U(t; U_n)\to U(t; U_0)$  in $H$
			for each $t>0$.
			\item[(iii)]
			The energy  equality
			\begin{align}\label{energy}
				&\sE(U(t))+\kappa\int_0^t \int_0^{L_0}\psi_{x\tau}^2 dx d\tau + \int_0^t \int_0^{L_0}\omega_{\tau}^2 dx d\tau \nonumber\\&=\sE(U_0)
				+\sum_{i=1, 2}\int_0^t   \int_0^{L_0} g_i  U_{it} dx d\tau +\sum_{i=3, 4}\int_0^t   \int_{L_0}^L g_i  U_{it} dx d\tau
			\end{align}
			holds
			for every $t>0$, where the energy functional $\sE$ is defined
			by the relation
			\begin{align*}
				\sE(U(t))=&\frac12\left[\int_0^{L_0}(\beta_1\phi_{t}^2+\rho_1\omega_{t}^2)dx+\int_{L_0}^L(\beta_2 u_{t}^2 +\rho_2 v_{t}^2)dx+\mu_2 \int_{L_0}^L u_{tx}^2 dx \right.\\&\left.+\mu_1 \int_0^{L_0} \phi_{tx}^2 dx+\lambda_1\int_0^{L_0} \phi_{xx}^2 dx+\lambda_2\int_{L_0}^L u_{xx}^2 dx+ Q_1(Z)+Q_2(Z)\right]
			\end{align*}
			with
			\begin{align}\label{Q-def}
				Q_1(Z)=\int_0^{L_0}\left(\omega_x+\frac{\phi_x^2}{2}\right)^2 dx,\;\;Q_2(Z)=\int_{L_0}^L\left(v_x+\frac{u_x^2}{2}\right)^2 dx.
			\end{align}
			\item[(vi)]
			If $g=0$ and $(Z_0, Z_1)\in D = (\EuScript D(G^{-1}N_1)\times \EuScript D(N_2))\times V$, weak solutions are strong, i.e. $Z(t)\in L_\infty(0,T;\EuScript D(G^{-1}N_1)\times \EuScript D(N_2))$, $Z_t(t)\in L_\infty(0,T;V)$, $Z_{tt}(t)\in L_\infty(0,T;W)$.
		\end{enumerate}
	\end{theorem}
	\begin{proof}
		The proof is quite standard, here we present the sketch of it. The existence of weak solutions can be shown by using the Galerkin  method  and relying on Lemma \ref{1}.
		We choose orthonormal bases $\tilde e^i=(e_1^i, e_2^i)$ in $X$ and $\hat e^i=(e_3^i, e_4^i)$ in $Y$
		consisting of eigenvectors of  operators $N_1$ and $N_2$.
		We define an approximate solution $Z_m=(\xi_m, \zeta_m)$, where  $\xi_m=(\phi_m,  u_m)=\sum_{i=1}^m d_i(t)\tilde e^i$ and            $\zeta_m=(\omega_m,  v_m)=\sum_{i=1}^m h_i(t)\hat e^i$ satisfying for $i=\overline {1, m}$
		\begin{align}
			\label{gal}
			&\beta_1 \int_0^{L_0}\phi_{mtt} e_1^i dx+ \mu_1\int_0^{L_0} \phi_{mttx} e_{1x}^i dx+\beta_2  \int_{L_0}^L u_{mtt} e_2^i dx+\mu_2 \int_{L_0}^L u_{mttx} e_{2x}^i dx \nonumber\\& \quad+\lambda_1 \int_0^{L_0}
			\phi_{mxx} e_{1xx}^i dx+\lambda_2
			\int_{L_0}^L u_{mxx} e_{2xx}^i dx+\kappa \int_0^{L_0} \phi_{mtx} e_{1x}^i dx \nonumber\\&\quad+ \int_0^{L_0}\left[\phi_{mx}\left(\omega_{mx}+1/2\phi_{mx}^2\right) \right] e_{1x}^i dx+ \int_{L_0}^L\left[u_{mx}\left(v_{mx}+1/2u_{mx}^2\right) \right] e_{2x}^i dx \nonumber\\&=\int_0^{L_0} g_1 e_1^i dx+\int_{L_0}^L g_2 e_2^i dx
		\end{align}
		and
		\begin{align}
			\label{gal1}
			&\rho_1 \int_0^{L_0}\omega_{mtt} e_3^i dx+\rho_2  \int_{L_0}^L v_{mtt} e_4^i dx+\gamma \int_0^{L_0}\omega_{mt} e_3^i dx+  \int_{L_0}^L\left(\omega_{mx}+1/2\phi_{mx}^2\right) e_{3x}^i dx \nonumber\\&\quad+\int_{L_0}^L\left(v_{mx}+1/2u_{mx}^2\right) e_{4x}^i dx=\int_0^{L_0} g_3 e_3^i dx+\int_{L_0}^L g_4 e_4^i dx
		\end{align}
		with initial conditions
		\[Z_m(0)=(\phi_m(0),  u_m(0), \omega_m(0), v_m(0))=(\phi_{0m}, u_{0m}, \omega_{0m}, v_{0m})=Z_{0m},\]
		\[Z_{tm}(0)=(\phi_{tm}(0), u_{tm}(0), \omega_{tm}(0), v_{tm}(0))=(\phi_{1m}, u_{1m}, \omega_{1m}, v_{1m})=Z_{1m}.\]
		Multilplying \eqref{gal} by $d_i'(t)$, \eqref{gal1} by $h_i'(t)$ and summing up with respect to $i$ from 1 to m we get
		\begin{align}\label{energy2}
			&\sE(Z_m(t), Z_{mt}(t))+\kappa \int_0^t \int_0^{L_0}\phi_{m x\tau}^2 dx d\tau+\gamma \int_0^t \int_0^{L_0}\omega_{m \tau}^2 dx d\tau   \nonumber\\&=\sE(Z_{0m}, Z_{1m})
			+\sum_{i=1, 2}\int_0^t   \int_0^{L_0} g_i  Z_{it} dx d\tau +\sum_{i=3, 4}\int_0^t   \int_{L_0}^L g_i  Z_{it} dx d\tau.
		\end{align}
		Using the Gronwall's lemma and Lemma 1 one can easily infer the estimate
		\begin{equation}\label{estgal}\|(Z_m(t),Z_{mt}(t))\|_{H}\le C(T,\|U_0\|_H,\|g\|_{W'})\end{equation}
		and the following convergences
		\begin{align}
			&\phi_m\to \phi, \;\text{weak-* in }\;L_\infty(0, T; H_{\{0\}}^2(0, L_0)), \label{sc1}\\
			&\omega_m\to \omega, \; \phi_{mt}\to \phi_t,\;  \text{weak-* in }\;L_\infty(0, T; H_{\{0\}}^1(0, L_0)), \label{sc2}\\
			&\omega_{mt}\to \omega_t,\;\;\text{weak-* in }\;L_\infty(0, T; L_2(0, L_0)),\label{sc3}\\
			&u_m\to u,\;\text{weak-* in}\;L_\infty(0, T; H^2\cap H_{\{L\}}^1(L_0, L)),\label{sc4}\\&v_m\to v,\; u_{tm}\to u_t,\;\;\text{weak-* in}\;\;L_\infty(0, T; H_{\{L\}}^1(L_0, L)), \label{sc5}\\ &v_{mt}\to v_t,\;\text{weak-* in }\;L_\infty(0, T; L_2(L_0, L)),\label{sc6}
		\end{align}
		and, consequently,  for any $\epsilon>1$
		\begin{align}
			&\phi_m\to \phi, \;\;\text{strongly in }\;L_2(0, T; H^{2-\epsilon}(0, L_0)),\label{cc7}\\
			&\omega_m\to \omega, \;\;\text{strongly in }\;L_2(0, T; H^{1-\epsilon}(0, L_0)),\label{cc8}\\
			&u_m\to u,\;\text{strongly in}\;L_2(0, T; H^{2-\epsilon}(L_0, L)),\label{cc9}\\
			&v_m\to v,\;\text{strongly in}\;L_2(0, T; H^{1-\epsilon}(L_0, L)).\label{cc10}
		\end{align}

		To prove the existence of weak solutions, one can resort to the standard limit procedure in \eqref{gal}. We only specify it for von Karman nonlinear terms.
		For any $\Psi=(\Psi_1, \Psi_2, \Psi_3, \Psi_4)\in
		\EuScript L_T^0$ it  follows from  \eqref{gal}, \eqref{gal1} that
		\begin{align}
			\label{1gal}
			&-\beta_1 \int_0^{T}\int_0^{L_0}\phi_{mt} \Psi_{1l} dx dt-\mu_1\int_0^{L_0} \phi_{mtx} \Psi_{1lx} dx dt-\beta_2  \int_0^{T}\int_{L_0}^L u_{mt} \Psi_{2l} dx dt \nonumber\\&\quad-\mu_2 \int_0^{T}\int_{L_0}^L u_{mtx} \Psi_{2lx} dx dt+\lambda_1 \int_0^{T}\int_0^{L_0}
			\phi_{mxx} \Psi_{1lxx} dx dt \nonumber\\ &\quad+\lambda_2
			\int_0^{T}\int_{L_0}^L u_{mxx} e_{2lxx} dx dt \nonumber\\&\quad+\kappa \int_0^{T}\int_0^{L_0} \phi_{mtx} \Psi_{1lx} dx dt+ \int_0^{T}\int_0^{L_0}\left[\phi_{mx}\left(\omega_{mx}+1/2\phi_{mx}^2\right) \right] \Psi_{1lx} dx dt \nonumber\\&\quad+ \int_0^{T}\int_{L_0}^L\left[u_{mx}\left(v_{mx}+1/2u_{mx}^2\right) \right] \Psi_{2lx} dx dt=\int_0^{T}\int_0^{L_0} g_1 \Psi_{1l} dx dt \nonumber\\&\quad+\int_0^{T}\int_{L_0}^L g_2 \Psi_{2l} dx dt+\beta_1 \int_0^{L_0}\phi_{mt} (0)\Psi_{1l}(0) dx-\mu_1\int_0^{L_0} \phi_{mtx}(0) \Psi_{1lx}(0) dx \nonumber\\&\quad-\beta_2  \int_{L_0}^L u_{mt}(0) \Psi_{2l}(0) dx-\mu_2 \int_{L_0}^L u_{mtx}(0)\Psi_{2lx}(0)dx
		\end{align}
		and
		\begin{align}
			\label{1gal1}
			&-\rho_1 \int_0^{T}\int_0^{L_0}\omega_{mt} \Psi_{3l} dx dt-\rho_2  \int_0^{T}\int_{L_0}^L v_{mt} \Psi_{4l} dx dt+\gamma \int_0^{T}\int_0^{L_0}\omega_{mt} \Psi_{3l} dx dt \nonumber\\&\quad+  \int_0^{T}\int_{L_0}^L\left(\omega_{mx}+1/2\phi_{mx}^2\right) \Psi_{3lx} dx dt+\int_0^{T}\int_{L_0}^L\left(v_{mx}+1/2u_{mx}^2\right) \Psi_{4lx} dx dt\nonumber\\&=\int_0^{L_0} g_3 \Psi_{3l} dx dt+\int_{L_0}^L g_4 \Psi_{4l} dx dt+\rho_1 \int_0^{L_0}\omega_{mt}(0) \Psi_{3l}(0) dx \nonumber\\& \quad+\rho_2  \int_{L_0}^L v_{mt}(0) \Psi_{4l}(0) dx,
		\end{align}
		where  $(\Psi_{1l}, \Psi_{2l})$ and $(\Psi_{3l}, \Psi_{4l})$ are  orthoprojections of $(\Psi_{1}, \Psi_{2})$ and $(\Psi_{3}, \Psi_{4})$ on the first basis vectors $\tilde e^i$ and $\hat e^i$ respectively, $l\le m$.
		Integrating by parts we get
		\begin{align*}
			&\int_{0}^T\int_0^{L_0}\left(\phi_{mx}\left(\omega_{mx}+1/2\phi_{mx}^2\right)-\phi_{x}\left(\omega_{x}+1/2\phi_{x}^2\right) \right)\Psi_{1lx} dx dt\\
			&\quad+\int_{0}^T\int_{L_0}^L\left(u_{mx}\left(v_{mx}+1/2u_{mx}^2\right)-u_{x}\left(v_{x}+1/2u_{x}^2\right) \right)\Psi_{2lx} dx dt\\ &=
			\int_{0}^T\int_0^{L_0}(\phi_{mx}\!-\!\phi_x)\omega_{x}\Psi_{1lx} dx dt\!-\!\int_{0}^T\int_0^{L_0}(\phi_{mxx}\Psi_{1lx}\!+\!\phi_{mx}\Psi_{1lxx})(\omega\!-\!\omega_m) dx dt\\&\quad+
			\int_{0}^T\int_{L_0}^L(u_{mx}\!-\!u_x)v_{x}\Psi_{2lx} dx dt\!-\!\int_{0}^T\int_0^{L_0}(u_{mxx}\Psi_{2lx}\!+\!u_{mx}\Psi_{2lxx})(v\!-\!v_m) dx dt\\&\quad+1/2\int_{0}^T\int_0^{L_0}(\phi_{mx}-\phi_x)(\phi_{mx}^2+\phi_{mx}\phi_x+\phi_x^2)\Psi_{1lx} dx dt\\&\quad+1/2\int_{0}^T\int_{L_0}^L(u_{mx}-u_x)(u_{mx}^2+u_{mx}u_x+u_x^2)\Psi_{2lx} dx dt
		\end{align*}
		For fixed $l$ we have from \eqref{cc7}--\eqref{cc10}
		\begin{align*}
			&\left| \int_{0}^T\int_0^{L_0}\left(\phi_{mx}\left(\omega_{mx}+1/2\phi_{mx}^2\right)-\phi_{x}\left(\omega_{x}+1/2\phi_{x}^2\right) \right)\Psi_{1lx} dx dt\right.\\&\quad\left.
			+\int_{0}^T\int_{L_0}^L\left(u_{mx}\left(v_{mx}+1/2u_{mx}^2\right)-u_{x}\left(v_{x}+1/2u_{x}^2\right) \right)\Psi_{2lx} dx dt\right|\\& \le
			C\left( \int_{0}^T\|\phi_{mx}-\phi_x\|\|\omega_{x}\|\|\Psi_{1lxx}\|  dt+\int_{0}^T\|\phi_{mxx}\|\|\Psi_{1lx}\|\|\omega-\omega_m\|  dt\right.\\& \quad+
			\int_{0}^T\|u_{mx}-u_x\|\|v_{x}\|\|\Psi_{2lxx}\| dt+\int_{0}^T\|u_{mxx}\|\|\Psi_{2lxx}\|\|v-v_m\| dt\\&\quad+\int_{0}^T\|\phi_{mx}-\phi_x\|(\|\phi_{mxx}\|^2+\|\phi_{xx}\|^2)\|\Psi_{1lxx} \| dt\\& \quad\left.+\int_{0}^T\|u_{mx}-u_x\|(\|u_{mxx}\|^2+\|u_{xx}\|^2)\|\Psi_{2lxx}\| dt\right)\to 0,\;\;m\to\infty.
		\end{align*}
		and
		\begin{align*}
			&|\int_{0}^T\int_0^{L_0} (\phi_{mx}^2-\phi_{x}^2) \Psi_{3lx} dx dt+\int_{0}^T\int_{L_0}^L (u_{mx}^2-u_{x}^2) \Psi_{4lx} dx dt |\\&\le  \int_{0}^T\|\phi_{mx}-\phi_{x}\|(\|\phi_{mxx}\|+\|\phi_{xx}\|) \|\Psi_{3lx}\| dt\\&\quad+\int_{0}^T\|u_{mx}-u_{x}\|(\|u_{mxx}\|+\|u_{xx}\|) \|\Psi_{4lx}\| dt  \to 0, \;\;m\to \infty.
		\end{align*}
		After the limit transition $l\to \infty$, we conclude that $Z=(\phi, u, \omega, v)$ is a weak solution to (\ref{1})-(\ref{24}). One can  prove the energy equality, uniqueness, continuity with respect to time and initial data using ideas presented in \cite{KochLa_2002}.\par
		Now we consider  the case $g=0$ and $Z_0\in \EuScript D(G^{-1}N_1)\times \EuScript D(N_2), \;Z_1\in V$. It is easy to see that $\tilde Z_m(t)=(\tilde \xi_m, \tilde \zeta_m)=(\tilde \phi_m, \tilde u_m, \tilde\omega_m, \tilde v_m)=Z_{mt}(t)$ satisfies
		\begin{align}
			\label{gals}
			&\beta_1 \int_0^{L_0}\tilde\phi_{mtt} e_1^i dx+ \mu_1\int_0^{L_0} \tilde\phi_{mttx} e_{1x}^i dx+\beta_2  \int_{L_0}^L \tilde u_{mtt} e_2^i dx+\mu_2 \int_{L_0}^L \tilde u_{mttx} e_{2x}^i dx \nonumber\\&\quad+\lambda_1 \int_0^{L_0}
			\tilde\phi_{mxx} e_{1xx}^i dx+\lambda_2
			\int_{L_0}^L \tilde u_{mxx} e_{2xx}^i dx+\kappa \int_0^{L_0} \tilde\phi_{mtx} e_{1x}^i dx\nonumber\\&\quad+ \int_0^{L_0}\left[\tilde \phi_{mx}\left(\omega_{mx}+1/2\phi_{mx}^2\right) \right] e_{1x}^i dx \nonumber\\&\quad+\int_{L_0}^L\left[\tilde u_{mx}\left(v_{mx}+1/2u_{mx}^2\right) \right] e_{2x}^i dx+
			\int_0^{L_0}\left[ \phi_{mx}\left(\tilde\omega_{mx}+\phi_{mx}\tilde\phi_{mx}\right) \right] e_{1x}^i dx \nonumber\\&\quad+\int_{L_0}^L\left[u_{mx}\left(\tilde v_{mx}+u_{mx}\tilde u_{mx}\right) \right] e_{2x}^i dx
			=0
		\end{align}
		and
		\begin{align}
			\label{gal1s}
			&\rho_1 \int_0^{L_0}\tilde\omega_{mtt} e_3^i dx\!+\!\rho_2  \int_{L_0}^L \tilde v_{mtt} e_4^i dx\!+\!\gamma \int_0^{L_0}\tilde \omega_{mt} e_3^i dx\!+\!  \int_{L_0}^L\left(\tilde\omega_{mx}+\tilde\phi_{mx}\phi_{mx}\right) e_{3x}^i dx \nonumber\\&\quad+\int_{L_0}^L\left(\tilde v_{mx}+\tilde u_{mx}u_{mx}\right) e_{4x}^i dx=0
		\end{align}
		with initial conditions
		\[\tilde Z_m(0)=(\tilde\phi_m(0),  \tilde u_m(0), \tilde\omega_m(0), \tilde v_m(0))=(\phi_{1m}, u_{1m}, \omega_{1m}, v_{1m})=Z_{1m}\in V,\]
		\[\tilde Z_{tm}(0)=(\tilde\phi_{tm}(0), \tilde u_{tm}(0), \tilde \omega_{tm}(0), \tilde v_{tm}(0))=(\phi_{2m}, u_{2m}, \omega_{2m}, v_{2m})\in W.\]
		Here
		\begin{align}
			&(\phi_{2m}, u_{2m})=-\tilde P_mG^{-1}(N_1(\phi_{m0}, u_{m0}) \nonumber\\&\qquad \qquad \qquad+([\phi_{0mx}(\omega_{0mx}+1/2\phi_{0mx}^2)]_x, [u_{0mx}(v_{0mx}+1/2u_{0mx}^2)]_x))\\
			& (\omega_{2m}, v_{2m})=-N_2(\omega_{m0}, v_{m0})+1/2\hat P_m([\phi_{0mx}^2]_x, [u_{0mx}^2]_x)),
		\end{align}
		where $\tilde P_m, \hat P_m$ are projectors on the first $m$ basis vectors $\tilde e_i$ and $\hat e_i$ respectively.
		It is easy to see that
		\begin{align}
			\label{1p}
			&\|G^{-1}N_1(\phi_{m0}, u_{m0})-G^{-1}N_1(\phi_{0}, u_{0})\|_{Y}^2 \nonumber\\&\le C (N_1((\phi_{m0}, u_{m0})-(\phi_{0}, u_{0}), G^{-1}N_1((\phi_{m0}, u_{m0})-(\phi_{0}, u_{0}))\nonumber\\&\le C \|N_1((\phi_{m0}, u_{m0})-(\phi_{0}, u_{0})\|_{Y'}^2 \nonumber\\&\le C\|N_1^{3/4}((\phi_{m0}, u_{m0})-(\phi_{0}, u_{0}))\|\to 0,\;m\to \infty
		\end{align}
		and
		\begin{align}
			\label{2p}
			&\|G^{-1}([\phi_{0mx}(\omega_{0mx}+1/2\phi_{0mx}^2)]_x, [u_{0mx}(v_{0mx}+1/2u_{0mx}^2)]_x \nonumber\\&\quad-[\phi_{0x}(\omega_{0x}+1/2\phi_{0x}^2)]_x, [u_{0x}(v_{0x}+1/2u_{0x}^2)]_x)\|_Y \nonumber\\&\le
			C(\|\phi_{0mxx}-\phi_{0xx}\|(\|\omega_{0xx}\|+\|\phi_{0xx}\|^2)+\|\phi_{0xx}\|\|\omega_{0mxx}-\omega_{0xx}\| \\& \quad+\|u_{0mxx}\!-\!u_{0xx}\|(\|v_{0xx}\|\!+\!\|u_{0xx}\|^2)\!+\!\|u_{0xx}\|\|v_{0mxx}\!-\!v_{0xx}\|)\!\to\! 0,\;m\!\to\! \infty.\nonumber
		\end{align}
		Analogously,
		\begin{align}
			\label{3p}
			&\|N_2(\omega_{m0}, v_{m0})-N_2(\omega_{0}, v_{0})+1/2([\omega_{0mx}^2]_x, [u_{0mx}^2]_x)-1/2([\omega_{0x}^2]_x, [u_{0x}^2]_x))\| \nonumber\\&\le
			C(\|v_{0mxx}-v_{0xx}\|+\|\omega_{0mxx}-\omega_{0xx}\|+\|\phi_{0mxx}-\phi_{0xx}\|\|\phi_{0xx}\| \nonumber\\&\quad+\|u_{0mxx}-u_{0xx}\|\|u_{0xx}\|))\to 0,\;\;m\to \infty.
		\end{align}
		Making use of \eqref{estgal}, \eqref{1p}--\eqref{3p}  it is easy to infer from \eqref{gals}--\eqref{gal1s} that
		\begin{equation}
			\label{1difgal}
			\|(\tilde Z_m(t), \tilde Z_{mt}(t))\|_{H}\le C(T,\|U_0\|_D).
		\end{equation}
		Then, it follows from \eqref{gal}--\eqref{gal1} that
		\begin{equation}
			\label{2difgal}
			\|Z_m(t)\|_{D}\le C(T,\|U_0\|_D).
		\end{equation}
		Arguing as above for weak solutions one can show the existence of strong solutions which are also approximate solutions by the energy argument.
	\end{proof}
	\begin{remark}
		The variational relation  in \eqref{sol_def} can be extended on the class of test functions from $\EuScript L_T$ by the appropriate limit transition, therefore, it is easy to see that in case $g_i=0$ strong solutions satisfy
		\begin{align}
			\label{48}
			&\int_0^T(\beta_1 \int_0^{L_0}\phi_{tt} z_1 dx+ \mu_1\int_0^{L_0} \phi_{ttx} z_{1x} dx+\beta_2 \int_{L_0}^L u_{tt} z_2 dx+\mu_2 \int_{L_0}^L u_{ttx} z_{2x} dx \nonumber\\&\quad-\lambda_1 \int_0^{L_0}
			\phi_{xxx} z_{1x} dx-\lambda_2
			\int_{L_0}^L u_{xxx} z_{2x} dx+\kappa \int_0^{L_0} \phi_{tx} z_{1x} dx \nonumber\\& \quad+ \int_0^{L_0}\left[\phi_{x}\left(\omega_{x}+1/2\phi_{x}^2\right) \right] z_{1x} dx \nonumber\\&\quad+\int_{L_0}^L\left[u_{x}\left(v_{x}+1/2u_{x}^2\right) \right] z_{2x} dx+
			\rho_1 \int_0^{L_0}\omega_{tt} y_1 dx+\gamma \int_0^{L_0}\omega_{t}y_1 dx \\&\quad+ \int_0^{L_0} \left(\omega_{x}+1/2\phi_{x}^2\right) y_{1x} dx+
			\rho_2 \int_{L_0}^L v_{tt} y_2dx+\int_{L_0}^L \left(v_{x}+1/2u_{x}^2\right) y_{2x} dx)dt=0  \nonumber
		\end{align}
		for any $(z_1, z_2),  (y_1, y_2)\in L_2(0,T;Y)$.
	\end{remark}
	\begin{remark}
		It follows from Theorem 3 that system \eqref{1}-\eqref{24} generates a dynamical system $(S_t, H)$ with the nonlinear operator $S_tU_0=U(t)$, where $U(t)$ is a weak solution to \eqref{1}-\eqref{24}, if $g_i=g_i(x), i=\overline{1, 4}$.
	\end{remark}
	\section{Asymptotic smoothness}
	To describe the long-time behaviour of solutions to system \eqref{1}--\eqref{24} we  show the following result on asymptotic smoothness.
	\begin{theorem}
		Let assumptions \eqref{g} hold true and  all right-hand sides in equations \eqref{1}--\eqref{4} be autonomous, i.e.
		\begin{equation}
			\label{auton}
			g_i=g_i(x), i=\overline{1, 4}.
		\end{equation}
		Let, moreover,
		\begin{equation}
			\label{koef}
			\beta_1\ge \beta_2,\;\rho_1\ge \rho_2,\;\mu_1\ge \mu_2, \;\lambda_1\le \lambda_2.
		\end{equation}
		Then the dynamical system $(S_t, H)$ generated by problem \eqref{1}--\eqref{24} is asymptotically smooth.
	\end{theorem}
	\begin{proof}
		Let $\tilde Z(t)=(\tilde\phi(t), \tilde u(t), \tilde\omega(t), \tilde v(t))$ and $\hat Z(t)=(\hat\phi(t), \hat u(t), \hat\omega(t), \hat v(t))$ be two weak solutions to problem \eqref{1}--\eqref{24} with initial data $\tilde U_0=(\tilde \phi_0, \tilde u_0, \tilde\omega_0, \tilde v_0, \tilde\phi_1, \tilde u_1,\\ \tilde\omega_1, \tilde v_1)$ and $\hat U_0=(\hat \phi_0, \hat u_0, \hat\omega_0,  \hat v_0, \hat\phi_1, \hat u_1, \hat\omega_1,  \hat v_1)$ respectively.  We assume that  $\tilde U_0$ and $\hat U_0$ lie in a ball $B_R$ of radius $R>0$. Then, it is easy to see from energy equality \eqref{energy} and Lemma \ref{lem1} that
		\begin{equation}
			\label{boun}
			\|(\tilde Z, \tilde Z_{t})\|_H+\|(\hat Z, \hat Z_{t})\|_H\le e^{cT}C(R).
		\end{equation}
		We consider the difference  $Z(t)\!=\!\tilde Z(t)\!-\!\hat Z(t)\!=\!(\phi(t), u(t), \omega(t), v(t))$ which satisfies  the  problem
		\begin{align*}
			&\beta_1 \phi_{tt}-\mu_1 \phi_{ttxx}-\kappa \phi_{txx}+\lambda_1
			\phi_{xxxx}\\&\qquad-\frac{1}{2}\left(\left[(\tilde \phi_x+\hat \phi_x)(\omega_x+(\tilde \phi_x+\hat \phi_x)\frac{\phi_x}{2})+\phi_x(\hat \omega_x+\tilde \omega_x+\frac{\hat \phi_x^2}{2}+\frac{\tilde \phi_x^2}{2}) \right]_x \right)=0,\\
			&\rho_1 \omega_{tt}-\gamma \omega_t-\left(\omega_x+\frac {\phi_x}{2}(\tilde \phi_x+\hat \phi_x)\right)_x=0,\\
			&\beta_2 u_{tt}-\mu_2 u_{ttxx}+\lambda_2
			u_{xxxx}\\&\qquad-\frac{1}{2}\left(\left[(\tilde u_x+\hat u_x)(v_x+(\tilde u_x+\hat u_x)\frac{u_x}{2})+u_x(\hat v_x+\tilde v_x+\frac{\hat u_x^2}{2}+\frac{\tilde u_x^2}{2}) \right]_x \right)=0,\qquad\qquad\qquad
			\\
			&\rho_2 v_{tt}-\left(v_x+\frac{u_x}{2}(\tilde u_x+\hat u_x)\right)_x=0,
			\\
			&\phi_{x}(L, t)=0, \; \phi(L, t)=0,\; \omega(L, t)=0,\;
			u_{xx}(L, t)=0, \; u(L, t)=0,\; v(L, t)=0,\\
			&\phi(L_0, t)=u(L_0, t),\;\omega(L_0, t)=v(L_0, t),\; \phi_x(L_0, t)=u_x(L_0, t),\;\\
			&\lambda_1 \phi_{xx}(L_0, t)
			=\lambda_2 u_{xx}(L_0, t),\\
			&\omega_x(L_0,t)=v_{x}(L_0,t),\;
			\left(\lambda_1\phi_{xxx}-\mu_1 \phi_{ttx}-\kappa \phi_{tx}\right)(L_0,t)
			=\left(\lambda_2 u_{xxx}-\mu_2 u_{ttx}\right)(L_0,t)
		\end{align*}
		with initial conditions $U_0=\tilde U_0-\hat U_0$ in a weak sense.
		
		First, by energy argument and integration over the interval $[t; T]$ we establish the following energy type equality
		\begin{equation}
			\label{en}
			\Phi(U(T))\!+\!\gamma\int_t^T \int_0^{L_0} \omega_{\tau}^2  dx d\tau\!+\! \kappa\int_t^T \int_0^{L_0} \phi_{x\tau}^2  dx d\tau \!=\! \Phi (U(t)) \!+\! \int_t^T H(\tilde Z, \hat Z) d\tau,
		\end{equation}
		where
		\begin{align*}
			\label{phi}
			\Phi(U(t))=&\frac12\bigg[\rho_1\int_0^{L_0}\omega_{t}^2 dx+
			\rho_2 \int_0^{L_0}v_{t}^2 dx+\beta_1 \int_{L_0}^L \phi_{t}^2dx+\beta_2 \int_{L_0}^L u_{t}^2dx\\&+\mu_1 \int_0^{L_0} \phi_{tx}^2 dx +\mu_2 \int_{L_0}^L u_{tx}^2 dx +\lambda_1\int_{L_0}^L \phi_{xx}^2 dx+ \int_{L_0}^L\omega_x^2 dx\\&\quad+\lambda_2\int_{L_0}^L u_{xx}^2 dx+\int_{L_0}^Lv_x^2 dx\bigg]
		\end{align*}
		and
		\begin{align}
			H(\tilde Z, \hat Z)&=\frac{1}{2}\int_0^{L_0} ((\tilde \omega_x+\hat \omega_x) \phi_x+\omega_x (\hat \phi_x+\tilde \phi_x)) \phi_{tx}dx+\int_0^{L_0} (\tilde \phi_x+\hat \phi_x) \phi_x \omega_{tx}dx \nonumber\\&\quad+\frac{1}{2}\int_0^{L_0}\!\! (\tilde  \phi_x^2\!+\!\tilde \phi_x \hat \phi_x\!+\!\hat \phi_x^2)\phi_x \phi_{tx}dx\!+\!
			\frac{1}{2}\int_{L_0}^L\!\! ((\tilde v_x\!+\!\hat v_x) u_x\!+\!v_x (\hat u_x\!+\!\tilde u_x)) u_{tx}dx \nonumber\\&\quad+\int_{L_0}^L (\tilde u_x+\hat u_x) u_x v_{tx}dx+\frac{1}{2}\int_{L_0}^L (\tilde  u_x^2+\tilde u_x \hat u_x+\hat u_x^2)u_x u_{tx}dx.
		\end{align}
		After integration of \eqref{en} over the interval $[0, T]$  we arrive at
		\begin{align}
			\label{stab2}
			&T\Phi(U(T)) +\gamma\int_0^T  \int_t^T \int_0^{L_0} \omega_{\tau}^2  dx d\tau dt  +\kappa\int_0^T  \int_t^T \int_0^{L_0} \phi_{x\tau}^2  dx d\tau dt \nonumber\\&= \int_0^T\Phi (U(t)) dt+ \int_0^T \int_t^T H(\tilde Z, \hat Z) d\tau dt.
		\end{align}
		The following computations can be justified by performing them on strong solutions and the limit procedure.
		We consider a function $\eta(x)\in C^\infty$ and assume that there exist  $0<\delta<L_0$, $\tilde\eta, \hat\eta>0$ such that
		\begin{align}
			\eta(0)=0, & \;\eta(L)=0, \label{prop1}\\
			-\tilde \eta\le\eta'(x)\le 0,& \; x\in (0, L_0-\delta),\\
			\eta'(x) >0, \;x\in (L_0-\delta, L),& \;\eta'(x)\ge \hat\eta >0, \; x\in (L_0-\frac{\delta}{2}, L).
		\end{align}
		We also choose  non-negative functions $\alpha(x), \sigma(x)\in C^\infty$  such that
		\begin{align}
			\alpha(L_0-\delta)= & 0 , \; \alpha'(L_0-\delta)=0,\\
			\alpha(x)=\hat\eta, &\; x\in (L_0-\frac{\delta}{2}, L),\\
			\alpha(x)=0, & \;x\in (0, L_0-\delta)
		\end{align}
		and
		\begin{align}
			\sigma(L_0)=0,&\; \sigma'(L_0)=0,\\
			\sigma(x)=\tilde \sigma=2\max &\{\hat\eta, \tilde\eta\}, x\in (0, L_0-\frac{\delta}{2}),\\
			\sigma(x)=0,& \;x\in (L_0, L).\label{prop2}
		\end{align}
		We choose  $(z_1, z_2)=(\eta(x) \phi_{x},\eta(x) u_{x})\in Y$  and  $(y_1, y_2)=(\eta(x) \omega_{x}, \eta(x) v_{x})\in Y$ in \eqref{48}, integrate over $[0, T]$, sum up the results and take the difference for solutions $\tilde Z(t)$, $\hat Z(t)$ and arrive at
		\begin{align}
			\label{63}
			&\beta_1 \int_0^T\int_0^{L_0} \phi_{tt}\eta(x)\phi_x dx dt+\mu_1 \int_0^T\int_0^{L_0} \phi_{ttx}\eta(x)\phi_{xx} dx dt \nonumber\\&\quad+\mu_1 \int_0^T\int_0^{L_0} \phi_{ttx}\eta'(x)\phi_{x} dx dt-\lambda_1 \int_0^T\int_0^{L_0} \phi_{xxx} \eta(x)\phi_{xx} dx dt\nonumber\\&\quad-\lambda_1 \int_0^T\int_0^{L_0} \phi_{xxx}\eta'(x)\phi_{x} dx dt+\rho_1 \int_0^T\int_0^{L_0} \omega_{tt}\eta(x)\omega_x dx dt \nonumber\\&\quad+\frac{1}{2}\int_0^T\int_0^{L_0} \left[(\tilde \phi_x+\hat \phi_x)(\omega_x+(\tilde \phi_x+\hat \phi_x)\frac{\phi_x}{2})\right] \eta(x)\phi_{xx} dx dt \nonumber\\&\quad+\frac{1}{2}\int_0^T\int_0^{L_0} \left[(\tilde \phi_x+\hat \phi_x)(\omega_x+(\tilde \phi_x+\hat \phi_x)\frac{\phi_x}{2})\right] \eta'(x)\phi_{x} dx dt \nonumber\\&\quad+
			\frac{1}{2}\int_0^T\int_0^{L_0}\left[\phi_x(\tilde \omega_x+\hat \omega_x+\frac{\tilde \phi_x^2}{2}+\frac{\hat \phi_x^2}{2})\right] \eta(x)\phi_{xx} dx dt \nonumber\\&\quad
			+\frac{1}{2}\int_0^T\int_0^{L_0}\left[\phi_x(\tilde \omega_x+\hat \omega_x+\frac{\tilde \phi_x^2}{2}+\frac{\hat \phi_x^2}{2})\right] \eta'(x)\phi_{x} dx dt \nonumber\\&\quad
			+\int_0^T\int_0^{L_0} (\omega_x+ (\tilde \phi_x+\hat \phi_x)\frac{\phi_x}{2})\eta(x)\omega_{xx} dx dt\nonumber\\&\quad+\int_0^T\int_0^{L_0} (\omega_x+ (\tilde \phi_x+\hat \phi_x)\frac{\phi_x}{2})_x \eta'(x)\omega_x dx dt
			\nonumber\\&\quad
			+
			\kappa \int_0^T\int_0^{L_0}\phi_{tx}\eta(x)\phi_{xx} dx dt+
			\kappa \int_0^T\int_0^{L_0}\phi_{tx}\eta'(x)\phi_{x} dx dt\nonumber\\&\quad+\gamma \int_0^T\int_0^
			{L_0}\omega_{t}\eta(x)\omega_x dx dt 
			+\beta_2 \int_0^T\int_{L_0}^L u_{tt}\eta(x)u_x dx dt \nonumber\\&\quad +\mu_2 \int_0^T\int_{L_0}^L u_{ttx}\eta(x)u_{xx} dx dt +\mu_2 \int_0^T\int_{L_0}^L u_{ttx}\eta'(x)u_{x} dx dt 
			\nonumber\\&\quad -\lambda_2 \int_0^T\int_{L_0}^L u_{xxx}\eta(x) u_{xx} dx dt-\lambda_2 \int_0^T\int_{L_0}^L u_{xxx} \eta'(x)u_{x} dx dt\nonumber\\&\quad+\rho_2 \int_0^T\!\int_{L_0}^L v_{tt}\eta(x)v_x dx dt +\frac{1}{2}\!\int_0^T\!\int_{L_0}^L \left[(\tilde u_x\!+\!\hat u_x)(v_x\!+\!(\tilde u_x\!+\!\hat u_x)\frac{u_x}{2})\right] \eta(x)u_{xx} dx dt \nonumber\\&\quad+\frac{1}{2}\int_0^T\int_{L_0}^L \left[(\tilde u_x+\hat u_x)(v_x+(\tilde u_x+\hat u_x)\frac{u_x}{2})\right] \eta'(x)u_{x} dx dt \nonumber\\&\quad+
			\frac{1}{2}\int_0^T\int_{L_0}^L\left[ u_x(\tilde v_x+\hat v_x+\frac{\tilde u_x^2}{2}+\frac{\hat u_x^2}{2})\right]\eta(x)u_{xx} dx dt \nonumber\\&\quad+
			\frac{1}{2}\int_0^T\int_{L_0}^L\left[ u_x(\tilde v_x+\hat v_x+\frac{\tilde u_x^2}{2}+\frac{\hat u_x^2}{2})\right]\eta'(x)u_{x} dx dt \nonumber\\&\quad+\int_0^T\int_{L_0}^L (v_x+(\tilde u_x+\hat u_x)\frac{u_x}{2}) \eta(x)v_{xx} dx dt\nonumber\\&\quad+\int_0^T\int_{L_0}^L (v_x+(\tilde u_x+\hat u_x)\frac{u_x}{2}) \eta'(x)v_x dx dt=0.
		\end{align}
		Integrating by parts we obtain
		\begin{align}
			\label{654}
			&\beta_1 \int_0^T\int_0^{L_0} \phi_{tt}\eta(x)\phi_x dx dt+\rho_1 \int_0^T\int_0^{L_0} \omega_{tt}\eta(x)\omega_x dx dt \nonumber\\&\quad+\beta_2 \int_0^T\int_{L_0}^L u_{tt}\eta(x)u_x dx dt
			+\rho_2 \int_0^T\int_{L_0}^L v_{tt}\eta(x)v_x dx dt \nonumber\\&=-\beta_1 \int_0^T\int_0^{L_0} \phi_{t}\eta(x)\phi_{xt}dx dt-\rho_1 \int_0^T\int_0^{L_0} \omega_{t}\eta(x)\omega_{tx} dx dt \nonumber\\&\quad+\beta_1 \int_0^{L_0} \phi_{t}(T)\eta(x)\phi_x(T)dx+\rho_1\int_0^{L_0} \omega_{t}(T)\eta(x)\omega_x(T)dx \nonumber\\&\quad-\beta_1 \int_0^{L_0} \phi_{t}(0)\eta(x)\phi_x(0)dx-\rho_1\int_0^{L_0} \omega_{t}(0)\eta(x)\omega_x(0)dx \nonumber\\&\quad-\beta_2 \int_0^T\int_{L_0}^L u_{t}\eta(x)u_{xt}dx dt-\rho_2 \int_0^T\int_{L_0}^L v_{t}\eta(x)v_{tx} dx dt \nonumber\\&\quad+\beta_2 \int_{L_0}^L u_{t}(T)\eta(x)u_x(T)dx+\rho_2\int_{L_0}^L v_{t}(T)\eta(x)v_x(T) dx\nonumber\\&\quad-\beta_2 \int_{L_0}^L u_{t}(0)\eta(x)u_x(0)dx 
			-\rho_2\int_{L_0}^L v_{t}(0)\eta(x)v_x(0) dx \nonumber\\&=\frac{\beta_1}{2}\int_0^T\int_0^{L_0} \phi_{t}^2 \eta'(x) dx dt+\frac{\rho_1}{2} \int_0^T\int_0^{L_0} \omega_{t}^2 \eta'(x)dx dt \nonumber\\&\quad
			+\frac{\beta_2}{2} \int_0^T\int_{L_0}^L u_{t}^2 \eta'(x)dx dt +\frac{\rho_2}{2} \int_0^T\int_{L_0}^L \eta'(x)v_{t}^2 dx dt\nonumber\\&\quad+(\beta_1-\beta_2)\frac{|\eta(L_0)|}{2} \int_0^T u_{t}^2(L_0) dt
			+(\rho_1-\rho_2)\frac{|\eta(L_0)|}{2} \int_0^T v_{t}^2(L_0) dt\nonumber\\&\quad+\beta_2 \int_{L_0}^L u_{t}(T)\eta(x)u_x(T)dx +\rho_2\int_{L_0}^L v_{t}(T)\eta(x)v_x(T) dx\nonumber\\&\quad-\beta_2 \int_{L_0}^L u_{t}(0)\eta(x)u_x(0)dx-\rho_2\int_{L_0}^L v_{t}(0)\eta(x)v_x(0) dx \nonumber\\&\quad\beta_1 \int_0^{L_0} \phi_{t}(T)\eta(x)\phi_x(T)dx+\rho_1\int_0^{L_0} \omega_{t}(T)\eta(x)\omega_x(T)dx\nonumber\\&\quad-\beta_1 \int_0^{L_0} \phi_{t}(0)\eta(x)\phi_x(0)dx-\rho_1\int_0^{L_0} \omega_{t}(0)\eta(x)\omega_x(0)dx.
		\end{align}
		Analogously, we have
		\begin{align}
			\label{641}
			&\mu_1 \int_0^T\int_0^{L_0} \phi_{ttx}\eta(x)\phi_{xx} dx dt+\mu_1 \int_0^T\int_0^{L_0} \phi_{ttx}\eta'(x)\phi_{x} dx dt \nonumber\\&\quad+\mu_2 \int_0^T\int_{L_0}^L u_{ttx}\eta(x)u_{xx} dx dt +\mu_2 \int_0^T\int_{L_0}^L u_{ttx}\eta'(x)u_{x} dx dt\nonumber\\& =
			-\frac{\mu_1}{2} \int_0^T\int_0^{L_0} \phi_{tx}^2\eta'(x)dx dt
			-\frac{\mu_2}{2} \int_0^T\int_{L_0}^L u_{tx}^2\eta'(x)dx dt \nonumber\\&\quad+(\mu_1-\mu_2)\frac{|\eta(L_0)|}{2} \int_0^T u_{tx}^2(L_0) dt\nonumber\\&\quad
			+\mu_1 \int_0^{L_0} \phi_{tx}(T)\eta'(x)\phi_{x}(T)dx +\mu_1 \int_0^{L_0} \phi_{tx}(T)\eta(x)\phi_{xx}(T)dx \nonumber\\&\quad
			-\mu_1 \int_0^{L_0} \phi_{tx}(0)\eta'(x)\phi_{x}(0)dx -\mu_1 \int_0^{L_0} \phi_{tx}(0)\eta(x)\phi_{xx}(0)dx \nonumber\\&\quad
			+\mu_2 \int_{L_0}^L u_{tx}(T)\eta'(x)u_{x}(T)dx +\mu_2 \int_{L_0}^L u_{tx}(T)\eta(x)u_{xx}(T)dx \nonumber\\&\quad -\mu_2 \int_{L_0}^L u_{tx}(0)\eta'(x)u_{x}(0)dx -\mu_2 \int_{L_0}^L u_{tx}(0)\eta(x)u_{xx}(0)dx.
		\end{align}
		It is easy to see that
		\begin{align}
			\label{164}
			&-\lambda_1 \int_0^T\int_0^{L_0} \phi_{xxx} \eta(x)\phi_{xx} dx dt-\lambda_1 \int_0^T\int_0^{L_0} \phi_{xxx}\eta'(x)\phi_{x} dx dt \nonumber\\&\quad-\lambda_2 \int_0^T\int_{L_0}^L u_{xxx} \eta(x) u_{xx} dx dt-\lambda_2 \int_0^T\int_{L_0}^L u_{xxx}\eta'(x) u_{x} dx dt \nonumber\\&
			=\frac{3\lambda_1}{2}\int_0^T\int_0^{L_0} \phi_{xx}^2 \eta'(x)dx dt+\frac{3\lambda_2}{2}\int_0^T\int_{L_0}^L u_{xx}^2\eta'(x) dx dt\nonumber\\&\quad+\frac{\lambda_2(\lambda_2-\lambda_1)|\eta(L_0)|}{2\lambda_1} \int_0^T u_{xx}^2(L_0) dt \nonumber\\&\quad+\lambda_1 \int_0^T\int_0^{L_0} \phi_{xx}\eta''(x)\phi_{x} dx+\lambda_2 \int_0^T\int_{L_0}^L u_{xx}\eta''(x) u_{x} dx dt
		\end{align}
		and
		\begin{align}
			\label{1644}
			&\frac{1}{2} \int_0^T\int_0^{L_0} \left[(\tilde \phi_x+\hat \phi_x)(\omega_x+ (\tilde \phi_x+\hat \phi_x)\frac{\phi_x}{2})\right] \eta(x)\phi_{xx} dx dt \nonumber\\&\quad+\int_0^T\int_0^{L_0} (\omega_x+(\tilde \phi_x+\hat \phi_x)\frac{\phi_x}{2}) \eta(x)\omega_{xx} dx dt\nonumber\\&\quad
			+\frac{1}{2}\int_0^T\int_{L_0}^L \left[(\tilde u_x+\hat u_x)(v_x+(\tilde u_x+\hat u_x)\frac{u_x}{2})\right] \eta(x) u_{xx} dx dt
			\nonumber\\&\quad+\int_0^T\int_{L_0}^L (v_x+(\tilde u_x+\hat u_x)\frac{u_x}{2}) \eta(x) v_{xx} dx dt \nonumber\\&
			=-\frac{1}{2} \int_0^T\int_0^{L_0} (\tilde \phi_{xx}+\hat \phi_{xx})(\omega_x+(\tilde \phi_x+\hat \phi_x)\frac{\phi_x}{2}) \eta(x)\phi_x dx dt \nonumber\\&\quad
			+\int_0^T\int_0^{L_0} (\omega_x+(\tilde \phi_x+\hat \phi_x)\frac{\phi_x}{2})_x \eta(x) (\omega_x+(\tilde \phi_x+\hat \phi_x)\frac{\phi_x}{2}) dx dt \nonumber\\&\quad-\frac{1}{2} \int_0^T\int_{L_0}^L (\tilde u_{xx}+\hat u_{xx})(v_x+(\tilde u_x+\hat u_x)\frac{u_x}{2}) \eta(x) u_x dx dt
			\nonumber\\&\quad+\int_0^T\int_{L_0}^L (v_x+(\tilde u_x+\hat u_x)\frac{u_x}{2})_x \eta(x) (v_x+(\tilde u_x+\hat u_x)\frac{u_x}{2}) dx dt \nonumber\\&
			=-\frac{1}{2} \int_0^T\int_0^{L_0} (\tilde \phi_{xx}+\hat \phi_{xx})(\omega_x+(\tilde \phi_x+\hat \phi_x)\frac{\phi_x}{2}) \eta(x)\phi_x dx dt
			\nonumber\\&\quad-\frac{1}{2}\int_0^T\int_0^{L_0}(\omega_x+(\tilde \phi_x+\hat \phi_x)\frac{\phi_x}{2})^2 \eta'(x) dx dt \nonumber\\&\quad-\frac{1}{2} \int_0^T\int_{L_0}^L (\tilde u_{xx}+\hat u_{xx})(v_x+(\tilde u_x+\hat u_x)\frac{u_x}{2}) \eta(x) u_x dx dt
			\nonumber\\&\quad-\frac{1}{2}\int_0^T\int_{L_0}^L (v_x+(\tilde u_x+\hat u_x)\frac{u_x}{2})^2 \eta'(x) dx dt .
		\end{align}
		Next we substitute   $(z_1, z_2)\!=\!(\alpha(x)\phi, \alpha(x) u)\!\in\!Y$,  $(y_1, y_2)\!=\!(0,0)$ into \eqref{48}, integrate over $[0, T]$, sum up the results and take the difference for solutions $\tilde Z(t)$, $\hat Z(t)$.
		\begin{align}\label{kl}
			&\beta_1 \int_0^T\int_0^{L_0}\phi_{tt}\alpha(x)\phi dx dt+\mu_1\int_0^T\int_0^{L_0} \phi_{ttx} \alpha(x)\phi_x dx dt\nonumber\\&\quad+\mu_1\int_0^T\int_0^{L_0} \phi_{ttx} \alpha'(x)\phi dx dt+\beta_2\int_0^T \int_{L_0}^L u_{tt} \alpha(x) u dx dt\nonumber\\&\quad+\mu_2 \int_0^T\int_{L_0}^L u_{ttx} \alpha(x) u_x dx dt+\mu_2 \int_0^T\int_{L_0}^L u_{ttx} \alpha'(x) u dx dt\nonumber\\&\quad-\lambda_1 \int_0^T \int_0^{L_0}
			\phi_{xxx} \alpha(x)\phi_x dx dt-\lambda_1 \int_0^T\int_0^{L_0} 
			\phi_{xxx} \alpha'(x)\phi dx dt \nonumber\\&\quad-\lambda_2
			\int_0^T\int_{L_0}^L u_{xxx} \alpha(x) u_x dx dt 
			-\lambda_2
			\int_0^T\int_{L_0}^L u_{xxx} \alpha'(x) u dx dt \nonumber
			\\&\quad+\kappa \int_0^T\int_0^{L_0} \phi_{tx}  \alpha(x)\phi_x dx dt+\kappa \int_0^T\int_0^{L_0} \phi_{tx}  \alpha'(x)\phi dx dt \nonumber\\&\quad
			+\frac{1}{2}\int_0^T\int_0^{L_0} \left[(\tilde \phi_x+\hat \phi_x)(\omega_x+(\tilde \phi_x+\hat \phi_x)\frac{\phi_x}{2})\right] \alpha(x)\phi_x dx dt \nonumber\\&\quad		
			+\frac{1}{2}\int_0^T\int_0^{L_0} \left[(\tilde \phi_x+\hat \phi_x)(\omega_x+(\tilde \phi_x+\hat \phi_x)\frac{\phi_x}{2})\right] \alpha'(x)\phi dx dt \nonumber\\&\quad
			+
			\frac{1}{2}\int_0^T\int_0^{L_0}\left[\phi_x(\tilde \omega_x+\hat \omega_x+\frac{\tilde \phi_x^2}{2}+\frac{\hat \phi_x^2}{2})\right] \alpha(x)\phi_x dx dt \nonumber\\&\quad  +
			\frac{1}{2}\int_0^T\int_0^{L_0}\left[\phi_x(\tilde \omega_x+\hat \omega_x+\frac{\tilde \phi_x^2}{2}+\frac{\hat \phi_x^2}{2})\right] \alpha'(x)\phi dx dt
			\nonumber\\&\quad
			+\frac{1}{2}\int_0^T\int_{L_0}^L \left[(\tilde u_x+\hat u_x)(v_x+(\tilde u_x+\hat u_x)\frac{u_x}{2})\right] \alpha(x)u_x dx dt\nonumber\\ &\quad+\frac{1}{2}\int_0^T\int_{L_0}^L \left[(\tilde u_x+\hat u_x)(v_x+(\tilde u_x+\hat u_x)\frac{u_x}{2})\right] \alpha'(x)u dx dt \nonumber\\&\quad+
			\frac{1}{2}\int_0^T\int_{L_0}^L\left[ u_x(\tilde v_x+\hat v_x+\frac{\tilde u_x^2}{2}+\frac{\hat u_x^2}{2})\right]\alpha(x)u_x dx dt \nonumber\\&\quad+\frac{1}{2}\int_0^T\int_{L_0}^L\left[ u_x(\tilde v_x+\hat v_x+\frac{\tilde u_x^2}{2}+\frac{\hat u_x^2}{2})\right]\alpha'(x)u dx dt=0.
		\end{align}
		Integrating by parts we come to
		\begin{align}
			\label{64}
			&\beta_1\int_0^T\int_0^{L_0} \alpha(x) \phi_{t}^2 dx dt+\beta_2\int_0^T\int_{L_0}^L \alpha(x) u_{t}^2 dx dt\nonumber\\&\quad+\mu_1 \int_0^T\int_0^{L_0} \alpha(x)  \phi_{tx}^2 dx dt+\mu_2 \int_0^T\int_{L_0}^L \alpha(x)  u_{tx}^2 dx dt \nonumber\\&\quad+\mu_1 \int_0^T\int_0^{L_0} \alpha'(x) \phi_{tx} \phi_{t}dx dt+
			\mu_2 \int_0^T\int_{L_0}^L \alpha'(x)  u_{tx} u_t dx dt \nonumber\\&\quad
			-\lambda_1\int_0^T\int_0^{L_0}\phi_{xx}^2\alpha(x) dx dt-\lambda_2\int_0^T\int_{L_0}^Lu_{xx}^2\alpha(x) dx dt \nonumber\\&\quad-\lambda_1\int_0^T\int_0^{L_0}\phi_{xx}\alpha'(x)\phi_x dx dt
			-\lambda_2\int_0^T\int_{L_0}^Lu_{xx}\alpha'(x) u_xdx dt\nonumber\\&\quad-\lambda_1\int_0^T\int_0^{L_0}\phi_{xx}\alpha''(x)\phi dx dt-\lambda_2\int_0^T\int_{L_0}^Lu_{xx}\alpha''(x) u dx dt\nonumber\\&\quad-\frac{1}{2}\int_0^T\int_0^{L_0} \left[(\tilde \phi_x+\hat \phi_x)(\omega_x+(\tilde \phi_x+\hat \phi_x)\frac{\phi_x}{2})\right] \alpha(x)\phi_x dx dt \nonumber\\&\quad-\frac{1}{2}\int_0^T\int_0^{L_0} \left[(\tilde \phi_x+\hat \phi_x)(\omega_x+(\tilde \phi_x+\hat \phi_x)\frac{\phi_x}{2})\right] \alpha'(x)\phi dx dt \nonumber\\&\quad
			-
			\frac{1}{2}\int_0^T\int_0^{L_0}\left[\phi_x(\tilde \omega_x+\hat \omega_x+\frac{\tilde \phi_x^2}{2}+\frac{\hat \phi_x^2}{2})\right] \alpha(x)\phi_x dx dt \nonumber\\&\quad-
			\frac{1}{2}\int_0^T\int_0^{L_0}\left[\phi_x(\tilde \omega_x+\hat \omega_x+\frac{\tilde \phi_x^2}{2}+\frac{\hat \phi_x^2}{2})\right] \alpha'(x)\phi dx dt \nonumber\\&\quad-\frac{1}{2}\int_0^T\int_{L_0}^L \left[(\tilde u_x+\hat u_x)(v_x+(\tilde u_x+\hat u_x)\frac{u_x}{2})\right] \alpha(x)u_x dx dt \nonumber\\&\quad-\frac{1}{2}\int_0^T\int_{L_0}^L \left[(\tilde u_x+\hat u_x)(v_x+(\tilde u_x+\hat u_x)\frac{u_x}{2})\right] \alpha'(x)u dx dt \nonumber\\&\quad
			-
			\frac{1}{2}\int_0^T\int_{L_0}^L\left[ u_x(\tilde v_x+\hat v_x+\frac{\tilde u_x^2}{2}+\frac{\hat u_x^2}{2})\right]\alpha(x)u_x dx dt \nonumber\\&\quad-
			\frac{1}{2}\int_0^T\int_{L_0}^L\left[ u_x(\tilde v_x+\hat v_x+\frac{\tilde u_x^2}{2}+\frac{\hat u_x^2}{2})\right]\alpha'(x)u dx dt\nonumber\\&
			=\beta_2  \int_{L_0}^L u_{t}(T) u(T) \alpha(x) dx+\beta_1  \int_0^{L_0} \phi_{t}(T) \phi(T) \alpha(x) dx \nonumber\\&\quad-\beta_2  \int_{L_0}^L u_{t}(0) u(0)\alpha(x) dx-\beta_1  \int_0^{L_0} u_{t}(0) u(0)\alpha(x) dx \nonumber\\&\quad-\mu_2 \int_{L_0}^L u_{tx}(T)u_{x}(T) \alpha(x)dx
			+\mu_2 \int_{L_0}^L u_{tx}(0) u_{x}(0) \alpha(x) dx \nonumber\\&\quad
			-\mu_1 \int_0^{L_0} \phi_{tx}(T)\phi_{x}(T) \alpha(x)dx+\mu_1 \int_0^{L_0} \phi_{tx}(0) \phi_{x}(0) \alpha(x) dx.
		\end{align}
		Next we substitute $(z_1, z_2)=(\sigma(x)\phi, 0)\in Y$ ,   $(y_1, y_2)=(\sigma(x)\omega, 0)\in Y$ into \eqref{48} , integrate over $[0, T]$, sum up the results, take the difference for solutions $\tilde Z(t)$, $\hat Z(t)$. After integration by parts we come to
		\begin{align}
			\label{66}
			&-\beta_1\int_0^T\int_0^{L_0} \sigma(x) \phi_{t}^2 dx dt-\mu_1 \int_0^T\int_0^{L_0} \sigma(x)  \phi_{tx}^2 dx dt \nonumber\\&\quad-\mu_1 \int_0^T\int_0^{L_0} \sigma'(x) \phi_{tx} \phi_{t}dx dt
			+\lambda_1\int_0^T\int_0^{L_0}\phi_{xx}^2\sigma(x) dx dt \nonumber\\&\quad+2\lambda_1\int_0^T\int_0^{L_0}\phi_{xx}\sigma'(x)\phi_x dx dt+\lambda_1\int_0^T\int_0^{L_0}\phi_{xx}\sigma''(x)\phi dx dt \nonumber\\&\quad+\frac{1}{2}\int_0^T\int_0^{L_0} \left[(\tilde \phi_x+\hat \phi_x)(\omega_x+(\tilde \phi_x+\hat \phi_x)\frac{\phi_x}{2})\right] \sigma(x)\phi_x dx dt \nonumber\\&\quad+\frac{1}{2}\int_0^T\int_0^{L_0} \left[(\tilde \phi_x+\hat \phi_x)(\omega_x+(\tilde \phi_x+\hat \phi_x)\frac{\phi_x}{2})\right] \sigma'(x)\phi dx dt \nonumber\\
			&\quad+
			\frac{1}{2}\int_0^T\int_0^{L_0}\left[\phi_x(\tilde \omega_x+\hat \omega_x+\frac{\tilde \phi_x^2}{2}+\frac{\hat \phi_x^2}{2})\right] \sigma(x)\phi_x dx dt \nonumber\\
			&\quad+
			\frac{1}{2}\int_0^T\int_0^{L_0}\left[\phi_x(\tilde \omega_x+\hat \omega_x+\frac{\tilde \phi_x^2}{2}+\frac{\hat \phi_x^2}{2})\right] \sigma'(x)\phi dx dt
			\nonumber\\&\quad
			-\rho_1\int_0^T\int_0^{L_0} \sigma(x)  \omega_{t}^2 dx dt-\gamma \int_0^T\int_0^{L_0} \sigma(x)  \omega_t\omega dx dt \nonumber\\&\quad+\int_0^T\int_0^{L_0} \sigma(x) \omega_x^2 dx dt+\frac{1}{2}\int_0^T\int_0^{L_0} \sigma(x) \omega_x\phi_x(\tilde \phi_x+\hat \phi_x)dx dt \nonumber\\&=\beta_1  \int_0^{L_0} \phi_{t}(0) \phi(0) \sigma(x) dx
			-\beta_1  \int_0^{L_0} \phi_{t}(T) \phi(T) \sigma(x) dx \nonumber\\&\quad
			+\mu_1 \int_0^{L_0} \phi_{tx}(T)\phi_{x}(T) \sigma(x)dx-\mu_1 \int_0^{L_0} \phi_{tx}(0) \phi_{x}(0) \sigma(x) dx \nonumber\\&\quad+\rho_1\int_0^{L_0} \sigma(x)\omega_{t}(T) \omega(T) dx-\rho_1\int_0^{L_0} \sigma(x)\omega_{t}(0) \omega(0) dx.
		\end{align}
		Collecting \eqref{koef}, \eqref{63}--\eqref{66} and taking into account properties of functions $\eta$, $\alpha$, $\sigma$ defined in \eqref{prop1}--\eqref{prop2} and using Lemma 2 and \eqref{boun} we come to
		\begin{multline}
			\label{659}
			\sum\limits_{i=1}^4\int_0^T\Psi_i(t) dt\le C(\Phi(0)+\Phi(T)+\int_0^T\int_0^{L_0} \phi_{xt}^2dxdt\\+\int_0^T\int_0^{L_0} \omega_{t}^2dxdt)+C(R,T)\int_0^T lot(\phi, \omega, u, v)dt,
		\end{multline}
		where
		\begin{align}
			\label{lot}
			lot(\phi, \omega, u, v)=&\|\phi\|_{H^{2-\epsilon}(0,L_0)}^2+\|u\|_{H^{2-\epsilon}(L_0,L)}^2 \nonumber\\&+\|\omega\|_{H^{1-\epsilon}(0,L_0)}^2+\|v\|_{H^{1-\epsilon}(L_0,L)}^2+\|u_t\|_{H^{1-\epsilon}(L_0,L)}^2
		\end{align}
		for $0<\epsilon<1/2$ and
		\begin{align*}
			\Psi_1(t)=&\hat\eta\bigg (\frac{3\beta_2}{2}\int_0^T\int_{L_0}^L u_{t}^2 dx dt +\frac{\mu_2}{2} \int_0^T\int_{L_0}^L u_{tx}^2 dx dt \nonumber\\&+\frac{\lambda_2}{4}\int_0^T\int_{L_0}^L u_{xx}^2dx dt
			+\frac{1}{4}\int_0^T\int_{L_0}^Lv_x^2 dx dt
			+\frac{\rho_2}{2} \int_0^T\int_{L_0}^L  v_{t}^2 dx dt \bigg),
			\\
			\Psi_2(t)=&\hat\eta\bigg (\frac{3\beta_2}{2}\int_0^T\int_{L_0-\delta/2}^{L_0} \phi_{t}^2 dx dt +\frac{\mu_1}{2} \int_0^T\int_{L_0-\delta/2}^{L_0} \phi_{tx}^2 dx dt \nonumber\\&+\frac{\lambda_1}{4}\!\!\int_0^T\!\!\!\int_{L_0-\delta/2}^{L_0}\!\! \phi_{xx}^2dx dt
			\!+\!\frac{1}{4}\!\int_0^T\!\!\!\int_{L_0-\delta/2}^{L_0}\omega_x^2 dx dt\!+\!\frac{\rho_1}{2} \!\int_0^T\!\!\!\int_{L_0-\delta/2}^{L_0}  \omega_{t}^2 dx dt \bigg),
			\\
			\Psi_3(t)=&\hat\eta\bigg (\frac{3\beta_2}{2}\int_0^T\int_{L_0-\delta}^{L_0-\delta/2} \phi_{t}^2 dx dt +\frac{\mu_1}{2} \int_0^T\int_{L_0-\delta}^{L_0-\delta/2} \phi_{tx}^2 dx dt \nonumber\\&+\frac{\lambda_1}{4}\!\int_0^T\!\!\int_{L_0-\delta}^{L_0-\delta/2}\!\! \phi_{xx}^2dx dt
			\!+\!\frac{1}{4}\!\int_0^T\!\!\int_{L_0-\delta}^{L_0-\delta/2}\!\!\omega_x^2 dx dt
			\!+\!\frac{\rho_1}{2} \!\int_0^T\!\!\int_{L_0-\delta}^{L_0-\delta/2} \!\! \omega_{t}^2 dx dt \bigg),
		\end{align*}
		and
		\begin{align*}
			\Psi_4(t)=&\tilde \eta\bigg (\frac{3\beta_2}{2}\int_0^T\int_0^{L_0-\delta} \phi_{t}^2 dx dt +\frac{\mu_1}{2} \int_0^T\int_0^{L_0-\delta} \phi_{tx}^2 dx dt \nonumber\\&\quad+\frac{\lambda_1}{4}\!\int_0^T\!\!\!\int_0^{L_0-\delta} \phi_{xx}^2dx dt\!+\!\frac{1}{4}\int_0^T\!\!\!\int_0^{L_0-\delta}\omega_x^2 dx dt
			\!+\!\frac{\rho_1}{2}\int_0^T\!\!\int_0^{L_0-\delta} \omega_{t}^2 dx dt \bigg).
		\end{align*}
		Consequently, \eqref{659} yields
		\begin{align}
			\label{6591}
			\int_0^T\Phi(t) dt\le& C(\Phi(0)+\Phi(T)+\int_0^T\int_0^{L_0} \phi_{xt}^2dxdt \nonumber\\&+\int_0^T\int_0^{L_0} \omega_{t}^2dxdt)+C(R,T)\int_0^T lot(\phi, \omega, u, v) dt,
		\end{align}
		where the last term is defined by \eqref{lot}.\par
		Then it follows from \eqref{en} with $t=0$ and \eqref{stab2}  that
		\begin{align}
			\label{70}
			\int_0^T\Phi(t) dt&\le C\left(\Phi(0)+\Phi(T)+\int_0^T\int_0^{L_0} \phi_{xt}^2dxdt+\int_0^T\int_0^{L_0} \omega_{t}^2dxdt\right) \nonumber\\& \quad+C(R,T)\int_0^T lot(\phi, \omega, \phi, u)dt \\&\le C\left(\Phi(0)+\Phi(T)+\left| \int_0^T H(\tilde Z, \hat Z)  d\tau\right|\right)+C(R, T) \int_0^Tlot(\phi, \omega, \phi, u) dt.\nonumber
		\end{align}
		Next we substitute \eqref{70} into \eqref{stab2}  and relying on \eqref{en} with $t=0$ come to the estimate
		\begin{align}
			\label{stab}
			&T\Phi(U(T)) +\int_0^T  \int_t^T \int_0^{L_0} \omega_{\tau}^2  dx d\tau dt+\int_0^T  \int_t^T \int_0^{L_0} \phi_{x\tau}^2  dx d\tau d \nonumber\\&\le C\left(\Phi(0)+\left|\int_0^T H(\tilde Z, \hat Z)  d\tau\right|+ \left|\int_0^T \int_t^T H(\tilde Z, \hat Z) d\tau dt\right|\right) \nonumber\\&\quad+C(R,T)\int_0^T lot(\phi, \omega, \phi, u)dt .
		\end{align}
		Now our remaining task is to estimate the nonlinear terms in \eqref{h}. We begin with the third term in the right-hand side.
		Integrating by parts with respect to $t$ and using Lemma 2 we obtain for any $0<\epsilon<1/2$
		\begin{align}
			\label{72}
			&\left|\int_0^T\int_{L_0}^L ((\tilde  u_x+\hat u_x)v_x u_{tx}+(\tilde u_x+\hat u_x) u_x v_{tx})dx dt\right| \nonumber\\&\le
			\left|\int_0^T\int_{L_0}^L  (\tilde  u_{xt}+\hat u_{xt})v_x u_{x}dx dt\right|+\max\limits_{[0,T]}\left|\int_{L_0}^L  (\tilde  u_{x}+\hat u_{x})v_x u_{x}dx \right| \nonumber
			\\&\le C(T)\max\limits_{[0,T]}( (\|\tilde  u_{xt}\|+\|\hat u_{xt}\|)\|v_x\|\|u\|_{2-\epsilon}+(\|\tilde  u_{xx}\|+\|\hat u_{xx}\|)\|v_x\|\|u\|_{2-\epsilon}) \nonumber\\&\le C(T, R)\max\limits_{[0,T]}\|u\|_{H^{2-\epsilon}(L_0,L)}.
			\\ \label{73}
			&\left|\int_0^T\int_{L_0}^L (\tilde v_x + \hat v_x)u_x u_{tx}dx dt\right|\le C
			\int_0^T(\|\tilde v_{x}\| +\| \hat v_{x}\|)\|u_{tx}\|\|u\|
			_{2-\epsilon}dt \nonumber\\&\le C(R, T)\max\limits_{[0,T]}\|u\|_{H^{2-\epsilon}(L_0,L)}.
		\end{align}
		Analogously,
		\begin{equation}\label{74}
			\left|\int_0^T\int_{L_0}^L (\tilde  u_x^2+\tilde u_x \hat u_x+\hat u_x^2)u_x u_{tx}dx\right|\le C(R, T)\max\limits_{[0,T]} \|u\|_{H^{2-\epsilon}(L_0,L)}.
		\end{equation}
		Collecting \eqref{72}--\eqref{74} and handling the first three terms in \eqref{h} in the same way and repeating the same computations for the first three terms  we infer from   \eqref{stab}
		\begin{equation}
			\label{stabb}
			\Phi(U(T)) \le \frac{C_1(R)}{T}+C_2(R, T) h(\phi, \omega, \phi, u),
		\end{equation}
		where
		\begin{align}
			\label{end}
			h(\phi, \omega, \phi, u)=& \max_{[0,T]}(\|\phi\|_{H^{2-\epsilon}(0,L_0)}^2+\|u\|_{H^{2-\epsilon}(L_0,L)}^2+\|\omega\|_{H^{1-\epsilon}(0,L_0)}^2 \\&+\|v\|_{H^{1\!-\!\epsilon}(L_0,L)}^2\!+\!\|u_t\|_{H^{1-\epsilon}(L_0,L)}^2\!+\!\|u\|_{H^{2-\epsilon}(L_0,L)}^2\!+\!\|\phi\|_{H^{2-\epsilon}(L_0, L)})\nonumber
		\end{align}
		for $0<\epsilon<1/2$.
		For any $\varepsilon>0$ it is possible to choose $T>0$ large enough to get \eqref{te} from \eqref{end} and to infer the statement of the theorem.
	\end{proof}
	\section {Existence of attractors}
	\begin{theorem} Let assumptions \eqref{g}, \eqref{auton}  hold true.
		Then, the dynamical system $(S_t, H)$ generated by \eqref{1}--\eqref{24} is gradient.
	\end{theorem}
	\begin{proof}
		
		{\it Step 1. Regularity.} Now we show that system \eqref{1}--\eqref{24} is gradient. The strict Lyapunov function is
		\begin{equation}
			\label{lap}
			\EuScript L(t)=\EuScript E(t)-\sum\limits_{i=1,2}\int_0^{L_0}g_iU_{i}dx-\sum\limits_{i=3, 4}\int_{L_0}^Lg_iU_{i}dx.
		\end{equation}
		In order to prove this, we assume that $\EuScript L(S_tU)=\EuScript L(U)$ , i.e. from the energy equality we obtain  $\phi_t=0$, $\omega_t=0$ and $\phi_d(t)=\phi(t+h)-\phi(t)=0$, $\omega_d(t)=\omega(t+h)-\omega(t)=0$ for arbitrary $T>h>0$.
		Consequently, $U_d(t)=(\phi_d(t), u_d(t),\omega_d(t), v_d(t))=U(t+h)-U(t)$, satisfies
		\begin{align}
			\label{sol_def_d}
			&-\int_0^T \int_{L_0}^L(\beta_2u_{dt} \Psi_{2t}+\rho_2v_{dt} \Psi_{4t})dx dt-\mu_2\int_0^T \int_{L_0}^L u_{dtx} \Psi_{2tx} dx dt\nonumber\\&\quad+\lambda_2 \int_0^T\int_{L_0}^{L}  u_{dxx} \Psi_{2xx}dx dt\nonumber\\&\quad+\int_0^T \int_{L_0}^L \left(v_{dx}\!+\!\frac{u_{dx}}{2}(u_x(t\!+\!h)\!+\!u_x(t)) \right)\left(\Psi_{4x}\!+\!\frac{1}{2}(u_x(t\!+\!h)\!+\!u_x(t)) \Psi_{2x}\right) dx dt \nonumber\\&\quad+\frac{1}{2}\int_0^T \int_{L_0}^L u_{dx}\left(v_x(t+h)+v_x(t)+\frac{u_x^2(t)}{2}+\frac{u_x^2(t+h)}{2} \right)\Psi_{2x}\ dx dt
			\nonumber\\&= \int_{L_0}^L(\beta_2u_{dt}(0) \Psi_{2}(0)+\rho_2v_{dt}(0) \Psi_{4}(0))dx+\mu_2 \int_{L_0}^L u_{dtx}(0) \Psi_{2x}(0) dx \nonumber\\&\quad
			- \int_{L_0}^L(\beta_2u_{dt}(T) \Psi_{2}(T)+\rho_2v_{dt}(T) \Psi_{4}(T))dx+\mu_2 \int_{L_0}^L u_{dtx}(T) \Psi_{2x}(T) dx.
		\end{align}
		Namely, $u_d$ and $v_d$ are weak solutions to an overdetermined problem
		\begin{align}
			\label{du}
			&\beta_2 u_{dtt}\!-\!\mu_2 u_{dttxx}+\lambda_2
			u_{dxxxx}\!-\!\frac{1}{2}\left[u_{dx}\left(v_x(t)\!+\!1/2u_x^2(t)\!+\!v_x(t+h)\!+\!1/2u_x^2(t\!+\!h)\right) \right]_x \nonumber\\&\quad-\frac{1}{2}\left[(u_{x}(t\!+\!h)\!+\!u_x(t))\left(v_{dx}\!+\!1/2u_{dx}(u_x(t\!+\!h)\!+\!u_x(t))\right) \right]_x\!=\!0,\; t\!>\!0,\;x\!\in\! (L_0, L),
			\\
			\label{dv}
			&\rho_2 v_{dtt}-\left(v_{dx}+1/2u_{dx}(u_x(t+h)+u_x(t))\right)_x=0
		\end{align}
		with boundary conditions
		\begin{align}
			&u_d(L_0, t)=0, \;v_d(L_0, t)=0,\; u_{dx}(L_0, t)=0,\;
			u_{dxx}(L_0, t)
			=0,\label{dd9}\\ v_{dx}(L_0,& t)=0,\;
			u_{dxxx} (L_0, t)= 0\;
			u_{dxx}(L, t)=0, \; u_d(L, t)=0,\; v_d(L, t)=0.\label{dd14}
		\end{align}
		Now we prove the additional regularity of the solution to the overdetermined problem \eqref{du}--\eqref{dd14}.
		We define the multiplier
		\begin{equation}
			\label{phii}
			\varphi(x,t)=(L-x)T^2-5(L-L_0)\left(t-\frac{T}{2}\right)^2.
		\end{equation}
		It is easy to see that
		\begin{equation}
			\label{phi1}
			\varphi(x, 0)\le -\frac{(L-L_0)T^2}{4}=-\sigma_0,\;\;\varphi(x, T)\le -\sigma_0.
		\end{equation}
		Moreover, there exist $t_0$ and $t_1$, such that and $0<t_0<\frac{T}{2}<t_1<T$
		\begin{equation}
			\label{phi2}
			\min\limits_{[t_0, t_1]}\varphi(x, t) \ge -\frac{\sigma_0}{2}.
		\end{equation}
		The following formal computations can be justified by performing them on the strong solutions to problem \eqref{du}, \eqref{dv} with boundary conditions
		\begin{align}
			&u_d(L_0, t)=0, \;v_d(L_0, t)=0,\; u_{dx}(L_0, t)=0,\nonumber\\
			&u_{dxx}(L, t)=0, \; u_d(L, t)=0,\; v_d(L, t)=0.\label{dd141}
		\end{align}
		and strong solutions to problem (\ref{1})-(\ref{24}) as the functions $u(t)$ and $v(t)$. It is easy to show by the energy arguments that weak solutions to \eqref{du}, \eqref{dv} can be approximated in the energy norm by such strong  solutions.
		For the simplification we present here the formal scheme and  substitute $\Psi_2=e^{\tau \varphi}(L-x)u_{dx}$ and $\Psi_4=e^{\tau \varphi}(L-x)v_{dx}$ into \eqref{sol_def_d}
		\begin{align}	\label{sol_def_d2}
			&-\int_0^T \int_{L_0}^L(\beta_2u_{dt}u_{dtx}+\rho_2v_{dt} v_{dtx})e^{\tau \varphi}(L-x)dx dt\nonumber\\&\quad+10\tau(L-L_0)\int_0^T \int_{L_0}^L(\beta_2u_{dt}u_{dx}+\rho_2v_{dt} v_{dx})(t-\frac{T}{2})e^{\tau \varphi}(L-x)dx dt\nonumber\\ &\quad-\mu_2\int_0^T \int_{L_0}^L u_{dtx} e^{\tau \varphi}(L-x)u_{dtxx} dx dt \nonumber\\&\quad+10\mu_2\tau(L-L_0)\int_0^T \int_{L_0}^L u_{dtx} e^{\tau \varphi}(t-\frac{T}{2})(L-x)u_{dxx} dx dt \nonumber\\&\quad-10\mu_2\tau(L-L_0)\int_0^T \int_{L_0}^L u_{dtx}(t-\frac{T}{2}) e^{\tau \varphi}[\tau T^2(L-x)+1]u_{dx} dx dt \nonumber\\&\quad+\mu_2\int_0^T \!\!\int_{L_0}^L u_{dtx}^2 e^{\tau \varphi}dx[\tau T^2(L\!-\!x)\!+\!1] dt
			\!+\!\lambda_2 \int_0^T\!\!\int_{L_0}^{L}  u_{dxx} u_{dxxx}e^{\tau \varphi}(L\!-\!x)dx dt\nonumber
			\\&\quad	
			-2\lambda_2 \int_0^T\int_{L_0}^{L}  u_{dxx}^2e^{\tau \varphi}[\tau T^2(L-x)+1]dx dt \nonumber\\&\quad+\lambda_2 \tau T^2\int_0^T\int_{L_0}^{L}  u_{dxx} u_{dx}e^{\tau \varphi}[\tau T^2(L-x)+2]dx dt \nonumber\\
			&\quad	+\int_0^T \int_{L_0}^L \left(v_{dx}+\frac{u_{dx}}{2}(u_x(t+h)+u_x(t)) \right)\nonumber\\&\qquad\left(v_{dx}+\frac{u_{dx}}{2}(u_x(t+h)+u_x(t)) \right)_x e^{\tau \varphi}(L-x)dx dt \nonumber\\&\quad-\int_0^T \int_{L_0}^L \left(v_{dx}+\frac{u_{dx}}{2}(u_x(t+h)+u_x(t)) \right)^2 [\tau T^2(L-x)+1]e^{\tau \varphi}dx dt \nonumber\\&\quad
			+\frac{1}{2}\int_0^T \int_{L_0}^L u_{dx}\left(v_x(t+h)+v_x(t)+\frac{u_x^2(t)}{2}+\frac{u_x^2(t+h)}{2} \right)u_{dxx}e^{\tau \varphi}(L-x) dx dt\nonumber\\&\quad-\frac{1}{2}\int_0^T\!\! \int_{L_0}^L u_{dx}^2\left(v_x(t\!+\!h)\!+\!v_x(t)\!+\!\frac{u_x^2(t)}{2}\!+\!\frac{u_x^2(t+h)}{2} \right)e^{\tau \varphi} [\tau T^2(L\!-\!x)\!+\!1] dx dt \nonumber\\&
			=\int_{L_0}^L(\beta_2u_{dt}(0) u_{dx}(0)+\rho_2v_{dt}(0)v_{dx}(0) )e^{\tau \varphi(0)}(L-x)dx \nonumber\\&\quad+\mu_2 \int_{L_0}^L u_{dtx}(0)u_{dxx}(0) e^{\tau \varphi(0)}(L-x) dx \nonumber\\&\quad
			-\mu_2 \int_{L_0}^L u_{dtx}(0)u_{dx}(0) e^{\tau \varphi(0)}[\tau T^2(L-x)+1] dx \nonumber\\&\quad-\int_{L_0}^L(\beta_2u_{dt}(T) u_{dx}(T)+\rho_2v_{dt}(T)v_{dx}(T) )e^{\tau \varphi(T)}(L-x)dx
			\nonumber\\&\quad-\mu_2 \int_{L_0}^L u_{dtx}(T)u_{dxx}(T) e^{\tau \varphi(T)}(L-x) dx
			\nonumber\\&\quad+\mu_2 \int_{L_0}^L u_{dtx}(T)u_{dx}(T) e^{\tau \varphi(T)}[\tau T^2(L-x)+1] dx,
		\end{align}
		After integration by parts we get
		\begin{align}
			\label{c2}
			&\int_0^T\!\!\! \int_{L_0}^L\!\! \left(v_{dx}\!+\!\frac{u_{dx}}{2}(u_x(t\!+\!h)\!+\!u_x(t)) \right)\left(v_{dx}\!+\!\frac{u_{dx}}{2}(u_x(t\!+\!h)\!+\!u_x(t)) \right)_x e^{\tau \varphi}(L\!-\!x)dx dt \nonumber\\&\quad-\int_0^T \!\!\!\int_{L_0}^L\!(\beta_2u_{dt}u_{dtx}\!+\!\rho_2v_{dt} v_{dtx})e^{\tau \varphi}\!(L\!-\!x)dx dt\!-\!\mu_2\!\!\int_0^T\!\!\!\int_{L_0}^L\!\! u_{dtxx}e^{\tau \varphi}\!(L\!-\!x) u_{dtx}dx dt \nonumber\\&\quad+\lambda_2 \int_0^T\int_{L_0}^{L}  u_{dxx} u_{dxxx}e^{\tau \varphi}(L-x)dx dt \nonumber\\&=-\frac{1}{2}\int_0^T \int_{L_0}^L(\beta_2u_{dt}^2+\rho_2v_{dt}^2)e^{\tau \varphi}[\tau T^2(L-x)+1]dx dt\nonumber\\&\quad-\frac{\mu_2}{2}\int_0^T\!\!\!\int_{L_0}^L\!\! u_{dtx}^2e^{\tau \varphi}[\tau T^2(L\!-\!x)\!+\!1]dx dt\!+\!\frac{\lambda_2}{2} \int_0^T\!\!\!\int_{L_0}^{L}\!\!  u_{dxx}^2e^{\tau \varphi}[\tau T^2(L\!-\!x)\!+\!1]dx dt\nonumber\\&\quad+
			\frac{1}{2}\int_0^T \int_{L_0}^L \left(v_{dx}+\frac{u_{dx}}{2}(u_x(t+h)+u_x(t)) \right)^2 e^{\tau \varphi}[\tau T^2(L-x)+1]dx dt \nonumber\\&\quad-\frac{\lambda_2}{2}(L-L_0)\int_0^Tu_{dxx}^2(L_0) e^{\tau \varphi(L_0)}dt-\frac{1}{2}(L-L_0)\int_0^Tv_{dx}^2(L_0) e^{\tau \varphi(L_0)}dt
		\end{align}
		Now we substitute  $(\Psi_2, \Psi_4)=(e^{\tau \varphi}[\tau T^2(L-x)+1]u_{d}, 0)$ into \eqref{du}
		\begin{align}
			\label{c111}
			&-\beta_2\int_0^T \int_{L_0}^Lu_{dt}^2 e^{\tau \varphi}[\tau T^2(L-x)+1]dx dt\nonumber\\&\quad+10\tau\beta_2(L-L_0)\int_0^T \int_{L_0}^Lu_{dt}u_d e^{\tau \varphi}[\tau T^2(L-x)+1](t-\frac{T}{2})dx dt\nonumber\\&\quad
			-\mu_2\int_0^T \int_{L_0}^L u_{dtx}^2e^{\tau \varphi} [\tau T^2(L-x)+1] dx dt\nonumber\\&\quad+\mu_2\tau T^2\int_0^T \int_{L_0}^L u_{dtx}u_{dt} e^{\tau \varphi}[\tau T^2(L-x)+2] dx dt \nonumber\\
			&\quad+10\mu_2\tau (L-L_0)\int_0^T \int_{L_0}^L u_{dtx}u_{dx} e^{\tau \varphi}[\tau T^2(L-x)+1](t-\frac{T}{2}) dx dt \nonumber\\&\quad-10\mu_2\tau^2 T^2 (L-L_0)\int_0^T \int_{L_0}^L u_{dtx}u_{d} e^{\tau \varphi}[\tau T^2(L-x)+2](t-\frac{T}{2}) dx dt \nonumber\\&\quad
			+\lambda_2 \int_0^T\int_{L_0}^{L}  u_{dxx}^2 e^{\tau \varphi}[\tau T^2(L-x)+1] dx dt\nonumber\\&\quad-2\lambda_2\tau T^2 \int_0^T\int_{L_0}^{L}  u_{dxx}  u_{dx} e^{\tau \varphi}[\tau T^2(L-x)+2] dx dt \nonumber\\
			&\quad+\lambda_2\tau^2 T^4 \int_0^T\int_{L_0}^{L}  u_{dxx}  u_{d} e^{\tau \varphi}[\tau T^2(L-x)+3] dx dt\nonumber\\&\quad+\int_0^T \int_{L_0}^L \left(v_{dx}+\frac{u_{dx}}{2}(u_x(t+h)+u_x(t)) \right)e^{\tau \varphi}\nonumber\\&\qquad \left(\frac{1}{2}(u_x(t\!+\!h)\!+\!u_x(t))\right)(u_{dx}[\tau T^2(L\!-\!x)\!+\!1]\!-\!u_d[\tau T^2(L\!-\!x)\!+\!2]) dx dt \nonumber\\
			&\quad+\frac{1}{2}\int_0^T\!\! \int_{L_0}^L u_{dx}^2\left(v_x(t\!+\!h)\!+\!v_x(t)\!+\!\frac{u_x^2(t)}{2}\!+\!\frac{u_x^2(t\!+\!h)}{2} \right) e^{\tau \varphi}[\tau T^2(L\!-\!x)\!+\!1]dx dt\nonumber\\&\quad-\frac{1}{2}\int_0^T\!\! \int_{L_0}^L u_{dx}u_{d}\left(v_x(t\!+\!h)\!+\!v_x(t)\!+\!\frac{u_x^2(t)}{2}\!+\!\frac{u_x^2(t\!+\!h)}{2} \right) e^{\tau \varphi}[\tau T^2(L\!-\!x)\!+\!2]dx dt
			\nonumber\\&= \beta_2\int_{L_0}^L(u_{dt}^2(0)e^{\tau \varphi(0)}-u_{dt}^2(T)e^{\tau \varphi(T)})[\tau T^2(L-x)+1]dx\nonumber\\&\quad+\mu_2 \int_{L_0}^L (u_{dtx}(0)u_{dx}(0)e^{\tau \varphi(0)}-u_{dtx}(T)u_{dx}(T)e^{\tau \varphi(T)})[\tau T^2(L-x)+1]dx \nonumber\\&\quad
			-\mu_2\tau T^2 \int_{L_0}^L (u_{dtx}(0)u_{d}(0)e^{\tau \varphi(0)}\!-\!u_{dtx}(T)u_{d}(T)e^{\tau \varphi(T)})[\tau T^2(L\!-\!x)\!+\!2]dx.
		\end{align}
		Distracting  \eqref{sol_def_d2} from \eqref{c111}  and taking into account \eqref{c2} we arrive at
		\begin{equation}
			\label{cm}
			L_d(T)+H_d(0)-H_d(T)=lot_d+M_d(T),
		\end{equation}
		where
		\begin{align}
			\label{ld}
			L_d(T)&=\frac{1}{2}\bigg(3\beta_2\int_0^T\int_{L_0}^L   e^{\tau \varphi}[\tau T^2(L-x)+1] u_{dt}^2 dx dt \nonumber\\&\quad+ \mu_2\int_0^T\int_{L_0}^L   e^{\tau \varphi}[\tau T^2(L-x)+1] u_{dtx}^2 dx dt \nonumber\\&\quad+\rho_2\int_0^T\!\!\int_{L_0}^L   e^{\tau \varphi}[\tau T^2(L\!-\!x)\!+\!1] v_{dt}^2 dx dt \nonumber\\&\quad+\!\lambda_2\int_0^T\!\!\int_{L_0}^L   e^{\tau \varphi}[\tau T^2(L\!-\!x)\!+\!1] u_{dxx}^2 dx dt\\ &\quad + \int_0^T\int_{L_0}^L   e^{\tau \varphi}[\tau T^2(L-x)+1] \left(v_{dx}+1/2u_{dx}(u_x(t+h)+u_x(t))\right)^2 dx dt \bigg),\nonumber
			\\
			\label{md}
			M_d(T)&=10 \mu_2\tau (L-L_0)\int_0^T\int_{L_0}^L(t-\frac{T}{2})u_{dtx}e^{\tau \varphi}(L-x) u_{dxx}dx dt\nonumber\\&\quad+10 \rho_2\tau (L-L_0)\int_0^T\int_{L_0}^L(t-\frac{T}{2})v_{dt}e^{\tau \varphi}(L-x) v_{dx}dx dt,
		\end{align}
		and
		\begin{align}
			\label{lotd}
			lot_d&=10\beta_2\tau(L-L_0)\int_0^T \int_{L_0}^Lu_{dt}u_{dx}(t-\frac{T}{2})e^{\tau \varphi}(L-x)dx dt\nonumber\\&\quad+\mu_2\tau T^2\int_0^T \int_{L_0}^L u_{dtx}u_{dt} e^{\tau \varphi}[\tau T^2(L-x)+2] dx dt\nonumber\\&\quad-10\mu_2\tau(L-L_0)\int_0^T \int_{L_0}^L u_{dtx} e^{\tau \varphi}[\tau T^2(L-x)+1](t-\frac{T}{2})u_{dx} dx dt\nonumber\\
			&\quad+\lambda_2 \tau T^2\int_0^T\int_{L_0}^{L}  u_{dxx} u_{dx}e^{\tau \varphi}[\tau T^2(L-x)+2]dx dt\nonumber\\&\quad
			-\int_0^T \int_{L_0}^L \left(v_{dx}+\frac{u_{dx}}{2}(u_x(t+h)+u_x(t)) \right) \nonumber\\&\qquad\left(v_{d}+\frac{u_{d}}{2}(u_x(t+h)+u_x(t)) \right) [\tau T^2(L-x)+1]e^{\tau \varphi}dx dt\nonumber\\&\quad
			+\frac{1}{2}\int_0^T \!\!\!\int_{L_0}^L\!\! u_{dx}\left(v_x(t\!+\!h)\!+\!v_x(t)\!+\!\frac{u_x^2(t)}{2}\!+\!\frac{u_x^2(t\!+\!h)}{2} \right)u_{dxx}e^{\tau \varphi}(L\!-\!x) dx dt \nonumber\\&\quad-\frac{1}{2}\int_0^T \!\!\!\int_{L_0}^L\!\! u_{dx}^2\left(v_x(t\!+\!h)\!+\!v_x(t)\!+\!\frac{u_x^2(t)}{2}\!+\!\frac{u_x^2(t+h)}{2} \right)e^{\tau \varphi} [\tau T^2(L\!-\!x)\!+\!1] dx dt \nonumber
			\\ &\quad
			+10\tau\beta_2(L-L_0)\int_0^T \int_{L_0}^Lu_{dt}u_d e^{\tau \varphi}[\tau T^2(L-x)+1](t-\frac{T}{2})dx dt \nonumber\\&\quad
			+10\mu_2\tau (L-L_0)\int_0^T \int_{L_0}^L u_{dtx}u_{dx} e^{\tau \varphi}[\tau T^2(L-x)+1](t-\frac{T}{2}) dx dt\nonumber\\&\quad-10\mu_2\tau^2 T^2 (L-L_0)\int_0^T \int_{L_0}^L u_{dtx}u_{d} e^{\tau \varphi}[\tau T^2(L-x)+2](t-\frac{T}{2}) dx dt \nonumber\\ &\quad
			-2\lambda_2\tau T^2 \int_0^T\int_{L_0}^{L}  u_{dxx}  u_{dx} e^{\tau \varphi}[\tau T^2(L-x)+2] dx dt \nonumber\\&\quad+\lambda_2\tau^2 T^4 \int_0^T\int_{L_0}^{L}  u_{dxx}  u_{d} e^{\tau \varphi}[\tau T^2(L-x)+3] dx dt \nonumber\\ &\quad+\int_0^T \int_{L_0}^L \left(v_{dx}+\frac{u_{dx}}{2}(u_x(t+h)+u_x(t)) \right)e^{\tau \varphi}\left(\frac{1}{2}(u_x(t+h)+u_x(t))\right) \nonumber\\&\qquad(u_{dx}[\tau T^2(L-x)+1]-u_d[\tau T^2(L-x)+2]) dx dt
			\\&\quad+\frac{1}{2}\int_0^T\!\!\! \int_{L_0}^L\!\! u_{dx}^2\left(v_x(t\!+\!h)\!+\!v_x(t)\!+\!\frac{u_x^2(t)}{2}\!+\!\frac{u_x^2(t\!+\!h)}{2} \right) e^{\tau \varphi}[\tau T^2(L\!-\!x)\!+\!1]dx dt \nonumber\\&\quad-\frac{1}{2}\int_0^T \!\!\!\int_{L_0}^L \!\!u_{dx}u_{d}\Big(v_x(t\!+\!h)\!+\!v_x(t)\!+\!\frac{u_x^2(t)}{2}\!+\!\frac{u_x^2(t\!+\!h)}{2} \Big) e^{\tau \varphi}[\tau T^2(L\!-\!x)\!+\!2]dx dt\nonumber
		\end{align}
		and
		\begin{multline}
			\label{hd}
			H_d(T)=\int_{L_0}^L(\beta_2u_{dt}(T) u_{dx}(T)+\rho_2v_{dt}(T)v_{dx}(T) )e^{\tau \varphi(T)}(L-x)dx\\+\mu_2 \int_{L_0}^L u_{dtx}(T)u_{dxx}(T) e^{\tau \varphi(T)}(L-x) dx\\+\beta_2\int_{L_0}^Lu_{dt}^2(T)e^{\tau \varphi(T)})[\tau T^2(L-x)+1]dx
			\\+\mu_2\tau T^2 \int_{L_0}^L u_{dtx}(T)u_{d}(T)e^{\tau \varphi(T)})[\tau T^2(L-x)+2]dx.
		\end{multline}
		It is easy to see that
		\begin{align}
			&10\mu_2\tau (L-L_0)\int_0^T\int_{L_0}^L(t-\frac{T}{2})u_{dtx}e^{\tau \varphi}(L-x) u_{dxx}dx dt \nonumber\\&\le
			\frac{\mu_2}{4}\int_0^T\int_{L_0}^L\tau T^2 u_{dtx}^2e^{\tau \varphi}(L-x) dx dt \nonumber\\&\quad+\frac{25\mu_2 (L-L_0)^2}{T^2}\int_0^T\int_{L_0}^L \tau T^2 e^{\tau \varphi}(L-x) u_{dxx}^2dx dt
		\end{align}
		and
		\begin{align}
			&10\rho_2\tau(L-L_0) \int_0^T\int_{L_0}^L v_{dt}e^{\tau \varphi} (L-x)(t-\frac{T}{2}) v_{dx} dx dt\\ &\le
			\frac{\rho_2}{4} \int_0^T\int_{L_0}^L \tau T^2 v_{dt}^2e^{\tau \varphi} (L-x) dx dt+\frac{25\rho_2(L-L_0)^2}{T^2} \int_0^T\int_{L_0}^L e^{\tau \varphi} (L-x) v_{dx}^2 dx dt.\nonumber
		\end{align}
		Taking into account \eqref{phi1} one can easily obtain that there exists $C_3>0$ such that
		\begin{equation}
			\label{hb}
			H_d(0)-H_d(T)\ge -C_3(E_d(T)+E_d(0))\tau^2e^{-\tau \sigma_0},
		\end{equation}
		where
		\begin{multline*}
			E_d(t)=\frac{1}{2}\left(\beta_2\int_{L_0}^L  u_{dt}^2 dx+ \mu_2\int_{L_0}^L  u_{dtx}^2 dx +\lambda_2\int_{L_0}^L  u_{dxx}^2 dx \right.\\\left.+ \int_{L_0}^L  \left(v_{dx}+1/2u_{dx}(u_x(t+h)+u_x(t))\right)^2 dx +\rho_2\int_{L_0}^L  v_{dt}^2 dx  \right).
		\end{multline*}
		Now we choose $T$ such that
		\begin{equation}\label{t}T>10 (L-L_0)\left(\frac{\mu_2}{\lambda_2}+\rho_2\right)^{1/2}.\end{equation}
		Consequently,  one can infer from \eqref{cm}--\eqref{t} the following estimate
		\begin{align}
			\label{cm1}
			&e^{-\frac{\tau \sigma_0}{2}}\int_{t_0}^{t_1}E_d(t)dt -C_3(E_d(T)+E_d(0))\tau^2e^{-\tau \sigma_0}\nonumber\\&\le C(\tau, T)\int_0^T\int_{L_0}^L (u_{dt}^2+u_{dx}^2+u_d^2+v_d^2)dx dt.
		\end{align}
		It is easy to infer from \eqref{sol_def_d} by the energy argument, Lemma 2, and the Gronwall's lemma that there exists $C(T, U_0)>0$ such that
		for $0\le s\le t\le T$
		\begin{equation}
			E_{d}(t)\le E_{d}(s)e^{C(T,U_0)(t-s)}\label{es11}
		\end{equation}
		and
		\begin{equation}
			E_{d}(s)\le E_{d}(t)e^{C(T,U_0)(t-s)}.\label{es22}
		\end{equation}
		If we choose $t=T$, $s=t$ in \eqref{es11} and $s=0$ in \eqref{es22}, we obtain
		\begin{equation}
			\label{p1}E_d(T)+E_d(0)\le CE_d(t)e^{C(T,U_0)}
		\end{equation}
		for any $t\in [0,T]$.
		
		Using \eqref{p1}  we estimate the left-hand side of
		\eqref{cm1} as follows
		\begin{align}
			\label{cm2}
			&e^{-\frac{\tau \sigma_0}{2}}\int_{t_0}^{t_1}E_d(t)dt -C_3(E_d(T)+E_d(0))\tau^2e^{-\tau \sigma_0}\nonumber\\&\ge \left(CT e^{-\frac{\tau \sigma_0}{2}}e^{-C(T,U_0)} -C_3 \tau^2e^{-\tau \sigma_0}\right)(E_d(T)+E_d(0)).
		\end{align}
		Choosing $\tau$ large enough and taking into account \eqref{cm1} we get
		\begin{equation}
			\label{cm11}
			E_d(T)+E_d(0)\le C(\tau, T)\int_0^T\int_{L_0}^L (u_{dt}^2+u_{dx}^2+u_d^2+v_d^2)dx dt.
		\end{equation}
		We note here that \eqref{cm11} was obtained for strong solutions to \eqref{sol_def_d} with boundary conditions \eqref{dd141}. Since these solutions are approximate in the energy norm for a weak solution for  overdetermined problem \eqref{cm11}, \eqref{dd9}, \eqref{dd14} one can get the same for the latter one by passing to the limit.\par
		Now we consider the solution $p\in H^4\cap H_0^1(L_0, L)$ to the elliptic problem
		\begin{align}
			\label{p}
			-p_{xx}&=u_d,\\
			p(L_0)=0,&\;\;p(L)=0.
		\end{align}
		Substituting $\Psi_2=p$, $\Psi_4=0$ into \eqref{sol_def_d} we get
		\begin{align}
			\label{oc}
			&\int_0^T \int_{L_0}^L u_{dt}^2 dx dt+\int_0^T \int_{L_0}^L p_{t}^2 dx dt \nonumber\\&\le C\bigg|\int_0^T\int_{L_0}^{L}  u_{dxx} u_{d}dx dt\nonumber\\&\quad
			+\int_0^T \int_{L_0}^L \Big(v_{dx}+\frac{u_{dx}}{2}(u_x(t+h)+u_x(t)) \Big)(u_x(t+h)+u_x(t))p_{x} dx dt\nonumber\\&\quad+\int_0^T \int_{L_0}^L u_{dx}\left(v_x(t\!+\!h)\!+\!v_x(t)\!+\!\frac{u_x^2(t)}{2}\!+\!\frac{u_x^2(t\!+\!h)}{2} \right)p_{x} dx dt\!+\! \int_{L_0}^Lu_{dt}(0) p(0))dx\nonumber\\&\quad+ \int_{L_0}^L u_{dtx}(0) p_{x}(0) dx
			+\int_{L_0}^Lu_{dt}(T) p(T))dx+ \int_{L_0}^L u_{dtx}(T) p_{x}(T) dx\bigg|\nonumber\\&\le C(\|U_0\|_H, T)\int_0^T\int_{L_0}^L (u_{dx}^2+v_d^2)dx dt +\varepsilon (E_d(0)+E_d(T))\nonumber\\&\quad+C(\varepsilon)(\|u_d(0)\|^2+\|u_d(T)\|^2).
		\end{align}
		Next after substituting \eqref{oc} into \eqref{cm11}, dividing  by $h^2$ and passing  to the limit $h\to 0$ we come to
		\begin{align}
			\label{cm111}
			E(U_t(T))&\le C( \|U_0\|_H, T)\left(\int_0^T\int_{L_0}^L (u_{tx}^2+u_t^2+v_t^2)dx dt+\|u_t(0)\|^2+\|u_t(T)\|^2\right)\nonumber\\&\le C( \|U_0\|_H, T).
		\end{align}
		Consequently, $(u_t,v_t, u_{tt}, v_{tt})\in H^2(L_0,L)\times H^1(L_0,L)\times H^1(L_0,L)\times L^2(L_0,L)$. Then, relying on \eqref{sol_def_d}, one can infer $(u_d, v_d)\in H^3(L_0, L)\times H^2(L_0, L)$ and, consequently, $(u_d, v_d, u_{dt}, v_{dt})\in H^3(L_0,L)\times H^2(L_0,L)\times H^2(L_0,L)\times H^1(L_0,L)$.\par
		Now we notice that $(u_t, v_t)=(w, z)$ is a weak solution to the problem
		\begin{align}
			\label{w}
			&\beta_2 w_{tt}-\mu_2 w_{ttxx}+\lambda_2
			w_{xxxx}-\left[w_x(v_x+1/2 u_x^2)\right]_x-\left[u_x\left(z_{x}+w_{x}u_x\right) \right]_x=0,\\
			\label{z}
			&\rho_2 z_{tt}-\left(z_{x}+w_{x}u_x\right)_x=0
		\end{align}
		with the same boundary conditions as in \eqref{dd9}--\eqref{dd14}.
		Consequently, $(w_d, z_d)=(w(t+h)-w(t), z(t+h)-z(t))$
		is a weak solution to the same problem with equations
		\begin{align}
			\label{wd}
			&\beta_2 w_{dtt}-\mu_2 w_{dttxx}+\lambda_2
			w_{dxxxx}-\left[w_{dx}(v_x(t+h)+1/2 u_x^2(t+h))\right]_x\nonumber\\&-\left[w_{x}(t)(v_{dx}\!+\!1/2 u_{dx}(u_x(t\!+\!h)\!+\!u_x(t)))\right]_x\!-\!\left[u_x(t)\left(z_{dx}\!+\!w_{dx}u_x(t)\!+\!w_{x}(t\!+\!h)u_{dx}\right) \right]_x\nonumber\\&\quad-\left[u_{dx}\left(z_{x}(x+h)+w_{x}(t+h)u_x(t+h)\right) \right]_x=0,
			\\
			\label{z1}
			&\rho_2 z_{dtt}-\left(z_{dx}+w_{dx}u_x(t)+ w_{x}(t+h)u_{dx}\right)_x=0.
		\end{align}
		Arguing as above  we come to
		\begin{equation}
			\label{cm22}
			E_{dd}(T)+E_{dd}(0)\le c(\tau, T)\int_0^T\int_{L_0}^L [w_{dx}^2+w_d^2+w_{dt}^2+z_d^2+v_{dx}^2+u_{dx}^2+u_{dxx}^2]dx dt,
		\end{equation}
		where
		\begin{multline*}
			E_{dd}(t)=\frac{1}{2}\left(\beta_2\int_{L_0}^L w_{dt}^2 dx+ \mu_2\int_{L_0}^L  w_{dtx}^2 dx +\lambda_2\int_{L_0}^L  w_{dxx}^2 dx \right.\\\left.+ \int_{L_0}^L  \left(z_{dx}+w_{dx}u_x(t))\right)^2 dx +\rho_2\int_{L_0}^L  z_{dt}^2 dx  \right).
		\end{multline*}
		Arguing as for \eqref{cm11} we get $(u_{tt},v_{tt}, u_{ttt}, v_{ttt})\in H^2(L_0,L)\times H^1(L_0,L)\times H^1(L_0,L)\times L^2(L_0,L)$ and, consequently, $(u_d,v_d, u_{dt}, v_{dt})\in H^4(L_0,L)\times H^3(L_0,L)\times H^3(L_0,L)\times H^2(L_0,L)$.\par
		{\it Step 2. Carleman estimates.}
		Now we derive  Carleman estimates (cf. \cite{isakov} for smoother solutions).   Let us consider the operator $P=\partial_{t}^2\partial_x -\rho\partial_x^3$, where $\rho$ is a positive constant.
		We choose $\tilde L>L$ and introduce functions
		\begin{equation}
			\label{eta1}
			r(x,t)=(x-\tilde L)^2 -m(t-\frac{T}{2})^2,\;\;\eta(x,t)=e^{\mu r(x,t)},
		\end{equation}
		\begin{equation}
			\label{eta2}
			q(x,t)=\tau \eta(x,t)\;\;\theta(x,t)=e^{q(x,t)},
		\end{equation}
		where $0\le t\le T$, $x\in (L_0, L )$. We select $T>0$ and $m\in (0, 1)$ as follows. Let
		\begin{equation}
			\label{T1}
			T>\max \{ \sqrt{4\tilde L^2+\sqrt{\rho}(\tilde L-L)}, \frac{4 \tilde L^2}{\sqrt{\rho}(\tilde L-L)}+1\}.
		\end{equation}
		Then,  there exists a constant $0<\sigma_1<\frac{\sqrt{\rho}(\tilde L-L)}{4}$ such that $T^2>4\tilde L^2+4\sigma_1$.
		Then we choose $m\in (0,1)$ such that
		\begin{equation}
			\label{T2}
			\frac{4\tilde L^2+4\sigma_1}{T^2}<m< \frac{\sqrt{\rho}(\tilde L-L)+4\sigma_1}{4T}<\frac{{\sqrt{\rho}(\tilde L-L)}}{2T}.
		\end{equation}
		This choice is possible due to \eqref{T1}.
		Therefore,
		\begin{equation}
			\label{q1}
			r(x,0)=r(x,T)=\tilde L^2-m\frac{T^2}{4}\le -\sigma_1
		\end{equation}
		for any $x\in (L_0, L)$.
		Moreover, there exist $t_0, t_1$  such that $0<t_0<\frac{T}{2}<t_1<T$ (chosen symmetrically around $\frac{T}{2}$)
		\begin{equation}
			\label{q2}
			\min\limits_{x\in [L_0, L], t\in [t_0, t_1]} r(x,t)\ge \sigma_2, \;\;0<\sigma_2<(\tilde L-L)^2.
		\end{equation}
		Now we set
		\begin{equation*}
			\tilde v=\theta v,
		\end{equation*}
		where we choose $v$ such that
		\begin{equation*}
			v\in H^3((L_0, L)\times (0, T)), \;v(L_0)=v_x(L_0)=v_{xx}(L_0)=v(L)=v_{xx}(L)=0.
		\end{equation*}
		Direct computations show that
		\begin{equation}
			\label{per}
			\theta Pv=\tilde v_{ttx}-\rho \tilde v_{xxx}+A_1 \tilde v_{xx}+A_2\tilde v_{x}+A_3\tilde v_{tt}+A_4\tilde v_{tx}+A_5\tilde v_t+A_6\tilde v,
		\end{equation}
		where
		\begin{align*}
			&A_1=3 \rho q_x,\;\;
			A_2=3\rho q_{xx}-3\rho (q_x)^2-q_{tt}+(q_t)^2,\\
			&A_3=-q_x,\;\;
			A_4=-2q_t,\;\;
			A_5=-2 q_{tx}+2 q_x q_t,\\
			&A_6=\rho q_{xxx}-3\rho q_x q_{xx} -q_{ttx}+2 q_t q_{tx}+q_x q_{tt}-q_x (q_t)^2+\rho (q_x)^3.
		\end{align*}
		We denote now
		\begin{align}
			&I_1=-\rho \tilde v_{xxx}+\tilde v_{ttx}+B_2 \tilde v_x+B_5 \tilde v_t,\label{1i}\\
			&I_2=B_1\tilde v_{xx}+B_4 \tilde v_{tx}+B_3  \tilde v_{tt}+B_6\tilde v,\label{2i}\\
			&S= S_2 \tilde v_x+S_5\tilde v_t+S_6\tilde v,\label{1s}
		\end{align}
		where
		\begin{align}
			&B_1=3 \rho q_x,\;\;
			B_2=-3\rho (q_x)^2+(q_t)^2,\;\;
			B_3=-q_x,\label{bb1}\\
			&B_4=-2q_t,\;\;
			B_5=2 q_x q_t,\;\;
			B_6=-q_x (q_t)^2+\rho (q_x)^3\label{bb2}
		\end{align}
		and
		\begin{align*}
			&S_2=3\rho q_{xx}-q_{tt},\;\;
			S_5=-2 q_{tx}\\
			&S_6=\rho q_{xxx}-3\rho q_x q_{xx}-q_{ttx}+2 q_t q_{tx}+q_x q_{tt}.
		\end{align*}
		Our next step is to estimate from below the expression
		\begin{align}
			\label{i2}
			&\int_0^T \int_{L_0}^{L} I_1 I_2 dx dt\nonumber\\&=-\rho\int_0^T \int_{L_0}^{L}  \tilde v_{xxx}B_1 \tilde v_{xx}dx dt-\rho\int_0^T \int_{L_0}^{L}  \tilde v_{xxx}B_4 \tilde v_{tx} dx dt\nonumber\\&\quad-\rho\int_0^T \int_{L_0}^{L}  \tilde v_{xxx}B_3 \tilde v_{tt} dx dt-\rho\int_0^T \int_{L_0}^{L}  \tilde v_{xxx}B_6 \tilde v dx dt+\int_0^T \int_{L_0}^{L}  \tilde v_{ttx}B_1 \tilde v_{xx} dx dt\nonumber\\&\quad+\int_0^T \int_{L_0}^{L}  \tilde v_{ttx}B_4 \tilde v_{tx} dx dt+\int_0^T \int_{L_0}^{L}  \tilde v_{ttx}B_3 \tilde v_{tt} dx dt+\int_0^T \int_{L_0}^{L}  \tilde v_{ttx}B_6 \tilde v dx dt\nonumber\\&\quad+\int_0^T \int_{L_0}^{L}  \tilde v_{x} B_1 B_2 \tilde v_{xx} dx dt+\int_0^T \int_{L_0}^{L}\tilde  v_{x} B_2 B_4 \tilde v_{tx} dx dt+\int_0^T \int_{L_0}^{L} \tilde v_{x} B_2 B_3 \tilde v_{tt}dx dt\nonumber\\&\quad+\int_0^T \int_{L_0}^{L}\tilde  v_{x} B_2 B_6 \tilde v dx dt+\int_0^T \int_{L_0}^{L} \tilde v_{t} B_5 B_1 \tilde v_{xx}dx dt+\int_0^T \int_{L_0}^{L}\tilde  v_{t} B_5 B_4 \tilde v_{tx} dx dt
			\nonumber\\&\quad+\int_0^T \int_{L_0}^{L} \tilde  v_{t} B_3 B_5 \tilde v_{tt} dx dt+\int_0^T \int_{L_0}^{L} \tilde  v_{t} B_5 B_6 \tilde v dx dt,
		\end{align}
		where $I_i$, $B_i$ are given by \eqref{1i}--\eqref{bb2}.
		Integrating by parts in the first three terms in the right-hand side of \eqref{i2} we obtain
		\begin{align}
			\label{b1}
			&-\rho\int_0^T \int_{L_0}^{L}  \tilde v_{xxx}B_1 \tilde v_{xx}dx dt\nonumber\\&=-3\rho^2 \int_0^T \int_{L_0}^{L}\tilde v_{xxx}  \tilde v_{xx} q_{x}dx dt=\frac{3}{2}\rho^2 \int_0^T \int_{L_0}^{L}  \tilde v_{xx}^2 q_{xx}dx dt,
			\\
			\label{b2}
			&-\rho\int_0^T \int_{L_0}^{L}  \tilde v_{xxx}B_4 \tilde v_{tx}dx dt\nonumber\\&=-2\rho\int_0^T \int_{L_0}^{L}  \tilde v_{xx} \tilde v_{txx}q_{t}dx dt-2\rho\int_0^T \int_{L_0}^{L}  \tilde v_{xx} \tilde v_{tx}q_{tx}dx dt\nonumber\\&=
			\rho\int_0^T \int_{L_0}^{L}  \tilde v_{xx}^2q_{tt}dx dt-\rho\int_{L_0}^{L}  \tilde v_{xx}^2(T)q_{t}(T)dx\nonumber\\&\quad+\rho\int_{L_0}^{L}  \tilde v_{xx}^2(0)q_{t}(0)dx-2\rho\int_0^T \int_{L_0}^{L}  \tilde v_{xx} \tilde v_{tx}q_{tx}dx dt,
		\end{align}
		and
		\begin{align}
			\label{b3}
			&-\rho\int_0^T \int_{L_0}^{L}  \tilde v_{xxx}B_3 \tilde v_{tt}dx dt\nonumber\\&=\rho\int_0^T \int_{L_0}^{L}\tilde v_{xxx} \tilde v_{tt} q_xdx dt\nonumber\\& = -\rho\int_0^T \int_{L_0}^{L}\tilde v_{xx} \tilde v_{ttx} q_xdx dt-\rho\int_0^T \int_{L_0}^{L}\tilde v_{xx} \tilde v_{tt} q_{xx} dx dt \nonumber\\&=-\frac{1}{2} \rho\int_0^T \int_{L_0}^{L} \tilde v_{tx}^2 q_{xx}dx dt+\rho\int_{L_0}^{L}  \tilde v_{xx}(0) \tilde v_{tx}(0)q_{x}(0)dx\nonumber\\&\quad-\rho\int_{L_0}^{L}  \tilde v_{xx}(T) \tilde v_{tx}(T)q_{x}(T)dx  -\rho\int_0^T \int_{L_0}^{L}\tilde v_{xx} \tilde v_{tt} q_{xx} dx dt \nonumber\\&\quad+\rho\int_0^T \int_{L_0}^{L}\tilde v_{xx} \tilde v_{tx} q_{tx} dx dt+\frac{1}{2} \rho\int_0^T  \tilde v_{tx}^2(L) q_{x}(L)dx dt.
		\end{align}
		Now we integrate by parts in the fourth term of the right-hand side of \eqref{i2}
		\begin{align}\label{b4}
			&-\rho\int_0^T \int_{L_0}^{L}  \tilde v_{xxx}B_6 \tilde v dx dt\nonumber\\&=-\rho\int_0^T \int_{L_0}^{L}\tilde v_{xxx} \tilde v (\rho q_x^3-q_x q_t^2 )dx dt\nonumber\\&=\rho\int_0^T \!\!\int_{L_0}^{L}\tilde v_{xx} \tilde v_x (\rho q_x^3\!-\!q_x q_t^2 )dx dt\!+\!
			\rho\int_0^T \!\!\int_{L_0}^{L}\tilde v_{xx} \tilde v (3\rho q_x^2 q_{xx}\!-\!q_{xx} q_t^2\!-\!2q_{x} q_t  q_{xt})dx dt
			\nonumber\\
			&=
			-\rho\int_0^T\!\! \int_{L_0}^{L}\tilde v_{x}^2  (\frac{9}{2}\rho q_x^2 q_{xx}\!-\!\frac{3}{2}q_{xx} q_t^2\!-\!3q_{x} q_t  q_{xt}) dx dt \nonumber\\&\quad+\!\frac{\rho}{2}\int_0^T \tilde v_{x}^2(L) q_x(L)(\rho q_x^2(L)\!-\!q_t^2(L) )dt\nonumber\\&\quad+
			\frac{\rho}{2}\int_0^T \!\!\int_{L_0}^{L}\tilde v^2 (3\rho q_x^2 q_{xxxx}\!+\!18\rho q_x q_{xx} q_{xxx}\!+\!6\rho q_{xx}^3\!-\!q_{xxxx} q_t^2\!-\!6 q_{xxx} q_t  q_{xt}\!-\!6 q_{xx} q_{tx}^2\nonumber\\&\quad-4q_{xx} q_t q_{txx}-6 q_x q_{tx} q_{txx}-2q_x q_t q_{txxx})dx dt.
		\end{align}
		Integrating by parts in the fifth, sixth, and seventh term in the right-hand side of \eqref{i2} we come to
		\begin{align}
			\label{b5}
			&\rho\int_0^T \int_{L_0}^{L}  \tilde v_{ttx}B_1 \tilde v_{xx} dx dt\nonumber\\&=3\rho\int_0^T \int_{L_0}^{L}  \tilde v_{ttx} \tilde v_{xx} q_x dx dt\nonumber\\&=-3\rho\int_0^T \int_{L_0}^{L}  \tilde v_{tx} \tilde v_{txx} q_x dx dt+3\rho \int_{L_0}^{L}  \tilde v_{tx}(T) \tilde v_{xx}(T) q_x(T) dx \nonumber\\&\quad-3\rho \int_{L_0}^{L}  \tilde v_{tx}(0) \tilde v_{xx}(0) q_x(0) dx  -3\rho\int_0^T \int_{L_0}^{L}  \tilde v_{tx} \tilde v_{xx} q_{xt} dx dt\nonumber\\&=
			\frac{3}{2}\rho\int_0^T \int_{L_0}^{L}  \tilde v_{tx}^2 q_{xx} dx dt-\frac{3}{2}\rho\int_0^T \tilde v_{tx}^2(L)q_x(L) dt\nonumber\\&\quad+3\rho \int_{L_0}^{L}  \tilde v_{tx}(T) \tilde v_{xx}(T) q_x(T) dx -3\rho \int_{L_0}^{L}  \tilde v_{tx}(0) \tilde v_{xx}(0) q_x(0) dx\nonumber\\&\quad-3\rho\int_0^T \int_{L_0}^{L}  \tilde v_{tx} \tilde v_{xx} q_{xt} dx dt,
			\\
			\label{b6}
			&\int_0^T \int_{L_0}^{L}  \tilde v_{ttx}B_4 \tilde v_{tx} dx dt\nonumber\\&=
			-2\int_0^T \int_{L_0}^{L}  \tilde v_{ttx} \tilde v_{tx} q_t dx dt\nonumber\\&=\int_0^T \int_{L_0}^{L}  \tilde v_{tx}^2  q_{tt} dx dt-\int_{L_0}^{L}   \tilde v_{tx}^2(T) q_t(T) dx + \int_{L_0}^{L}   \tilde v_{tx}^2(0) q_t(0) dx,
		\end{align}
		and
		\begin{equation}
			\label{b7}
			\int_0^T \int_{L_0}^{L}  \tilde v_{ttx}B_3 \tilde v_{tt} dx dt=-\int_0^T \int_{L_0}^{L}  \tilde v_{ttx} \tilde v_{tt} q_x dx dt=\frac{1}{2}\int_0^T \int_{L_0}^{L}  \tilde v_{tt}^2 q_{xx} dx dt.
		\end{equation}
		After integration by parts in the eighth term in the right-hand side of \eqref{i2} we get
		\begin{align}
			\label{b8}
			&\int_0^T \int_{L_0}^{L}  \tilde v_{ttx}B_6 \tilde v dx dt\nonumber\\&=\int_0^T \int_{L_0}^{L}  \tilde v_{ttx} \tilde v (\rho q_x^3-q_x q_t^2) dx dt\nonumber\\&=-\int_0^T\!\! \int_{L_0}^{L}  \tilde v_{tx} \tilde v_t (\rho q_x^3-q_x q_t^2) dx dt-\int_0^T \!\!\int_{L_0}^{L}  \tilde v_{tx} \tilde v (3\rho q_x^2 q_{tx}-q_{tx} q_t^2-2q_{x} q_t q_{tt}) dx dt\nonumber\\&\quad+\int_{L_0}^{L} \!\! \tilde v_{tx}(T) \tilde v(T)q_x(T) (\rho q_x^2(T)\!-\!q_t^2(T)) dx \!-\!\int_{L_0}^{L} \!\! \tilde v_{tx}(0) \tilde v(0)q_x(0) (\rho q_x^2(0)\!-\!q_t^2(0)) dx\nonumber\\&=\frac{1}{2}\int_0^T \int_{L_0}^{L}  \tilde v_{t}^2(3\rho q_x^2 q_{xx}-q_{xx} q_t^2-2q_x q_t q_{tx}) dx dt\nonumber\\&\quad+
			\int_0^T \int_{L_0}^{L}  \tilde v_{t} \tilde v_x (3\rho q_x^2 q_{tx}-q_{tx} q_t^2-2q_{x} q_t q_{tt}) dx dt\nonumber\\&\quad-\int_{L_0}^{L}\!\!  \tilde v_{t}(T) \tilde v_x(T)q_x(T) (\rho q_x^2(T)\!-\!q_t^2(T)) dx \!+\!\int_{L_0}^{L}  \tilde v_{t}(0) \tilde v_x(0)q_x(0) (\rho q_x^2(0)\!-\!q_t^2(0)) dx\nonumber\\&\quad -\int_{L_0}^{L}  \tilde v_{x}(T) \tilde v(T)(3\rho q_x^2(T) q_{xx}(T)-q_{xx}(T)q_t^2(T)-2q_{x}(T)q_t(T)q_{tx}(T)) dx \nonumber\\&\quad
			+\int_{L_0}^{L}  \tilde v_{x}(0) \tilde v(0) (3\rho q_x^2(0) q_{xx}(0)-q_{xx}(0)q_t^2(0)-2q_{x}(0)q_t(0)q_{tx}(0)) dx \nonumber\\&\quad-\frac{1}{2}\int_0^T \int_{L_0}^{L}  \tilde v^2 (3\rho q_x^2 q_{ttxx}+6\rho q_x q_{xx}q_{ttx}+6\rho q_{xx} q_{tx}^2+12\rho q_x q_{tx} q_{txx}-q_{ttxx} q_t^2\nonumber\\&\quad-2q_t q_{tx} q_{ttx}-2q_{tx}^2q_{tt}-2q_{xt} q_t q_{ttx}-2q_{txx} q_{tt} q_t-2q_{x} q_t q_{tttx}-2 q_x q_{tx} q_{ttt}\nonumber\\&\quad-2q_{xx} q_t q_{ttt}-2 q_{xx} q_{tt}^2-2 q_{xt}q_t q_{ttx}-4q_x q_{tt}q_{ttx}-2q_{txx}q_t q_{tt}-2q_{tt}q_{xt}^2) dx dt.
		\end{align}
		Next we integrate by parts in the ninth, the tenth, and the eleventh terms in the right-hand side of \eqref{i2}
		\begin{align}
			\label{b9}
			&\int_0^T \int_{L_0}^{L}  \tilde v_{x}B_2 B_1 \tilde v_{xx} dx dt\nonumber\\&= \int_0^T \int_{L_0}^{L}  \tilde v_{x} \tilde v_{xx}(-9\rho^2 q_x^3+3\rho q_x q_t^2) dx dt\nonumber\\&=\int_0^T \int_{L_0}^{L}  \tilde v_{x}^2(\frac{27}{2}\rho^2 q_x^2 q_{xx}-\frac{3}{2}\rho q_{xx} q_t^2-3\rho q_x q_t q_{tx}) dx dt \nonumber\\&\quad-\frac{3\rho}{2}\int_0^T  \tilde v_{x}^2(L) q_x(L)(3\rho q_x^2(L)- q_t^2(L)) dt,
			\\
			\label{b10}
			&\int_0^T \int_{L_0}^{L}  \tilde v_{x}B_2 B_4 \tilde v_{tx} dx dt\nonumber\\&=-2\int_0^T \int_{L_0}^{L}  \tilde v_{x} \tilde v_{tx}(q_t^3-3\rho q_t q_x^2) dx dt \nonumber\\&=-\int_{L_0}^{L}  \tilde v_{x}^2(T)(q_t^3(T)-3\rho q_t(T) q_x^2(T)) dx +\int_{L_0}^{L}  \tilde v_{x}^2(0)(q_t^3(0)-3\rho q_t(0) q_x^2(0)) dx \nonumber\\&\quad+\int_0^T \int_{L_0}^{L}  \tilde v_{x}^2 (3q_t^2 q_{tt}-3\rho q_{tt} q_x^2-6\rho q_t q_x q_{tx}) dx dt,
		\end{align}
		and
		\begin{align}
			\label{b11}
			\int_0^T \int_{L_0}^{L}  \tilde v_{x}B_2 B_3 \tilde v_{tt} dx dt=& \int_0^T \int_{L_0}^{L}  \tilde v_{x} \tilde v_{tt} (3\rho q_x^3-q_x q_t^2)dx dt\nonumber\\=&\frac{1}{2}\int_0^T \int_{L_0}^{L} \tilde v_{t}^2 (9\rho q_x^2 q_{xx}-q_{xx} q_t^2-2q_x q_t q_{tx})dx dt \nonumber
			\\&-\int_0^T \int_{L_0}^{L}  \tilde v_{x} \tilde v_{t} (9\rho q_x^2 q_{xt}-q_{xt} q_t^2-2q_x q_t q_{tt})dx dt \nonumber\\&+\int_{L_0}^{L}  \tilde v_{x}(T) \tilde v_{t}(T)q_x(T) (3\rho q_x^2(T)-q_t^2(T))dx\nonumber\\&-\int_{L_0}^{L}  \tilde v_{x}(0) \tilde v_{t}(0)q_x(0) (3\rho q_x^2(0)-q_t^2(0))dx.
		\end{align}
		Next we integrate by parts in the twelfth term in the right-hand side of \eqref{i2}
		\begin{align}
			\label{b12}
			&\int_0^T \int_{L_0}^{L}  \tilde v_{x}B_2 B_6 \tilde v dx dt \nonumber\\&= \int_0^T \int_{L_0}^{L}  \tilde v_{x} \tilde v (4\rho q_x^3 q_t^2-3\rho^2 q_x^5-q_x q_t^4) dx dt\nonumber\\&=-\int_0^T \!\!\!\int_{L_0}^{L}\! \tilde v^2 (6\rho q_x^2 q_{xx} q_t^2\!+\!4\rho q_x^3 q_t q_{tx}\!-\!\frac{15}{2}\rho^2 q_x^4 q_{xx}\!-\!\frac{1}{2}q_{xx} q_t^4\!-\!2q_x q_t^3 q_{tx}) dx dt.
		\end{align}
		Integrating by parts in the thirteenth, fourteenth, and fifteenth terms in the right-hand side of \eqref{i2} we obtain
		\begin{align}
			\label{b13}
			&\int_0^T \int_{L_0}^{L}  \tilde v_{t}B_1 B_5 \tilde v_{xx} dx dt\nonumber\\&=6\rho \int_0^T \int_{L_0}^{L}  \tilde v_{t} \tilde v_{xx} q_x^2 q_t dx dt\nonumber\\&=-6\rho \int_0^T \int_{L_0}^{L}  \tilde v_{tx} \tilde v_{x} q_x^2 q_t dx dt-6\rho \int_0^T \int_{L_0}^{L}  \tilde v_{t} \tilde v_{x}( q_x^2 q_{tx}+2q_x q_{xx} q_t )dx dt \nonumber\\
			&=3\rho \int_0^T \int_{L_0}^{L}  \tilde v_{x}^2 (2q_{tx}q_x q_t +q_x^2 q_{tt})dx dt+6\rho \int_0^T \int_{L_0}^{L}  \tilde v_{t} \tilde v_{x}( q_x^2 q_{tx}+2q_x q_{xx} q_t )dx dt\nonumber\\&\quad-3\rho \int_{L_0}^{L}  \tilde v_{x}^2(T) q_x^2(T) q_t (T)dx+3\rho \int_{L_0}^{L}  \tilde v_{x}^2(0) q_x^2(0) q_t (0)dx,
			\\
			\label{b14}
			&\int_0^T \int_{L_0}^{L}  \tilde v_{t}B_4 B_5 \tilde v_{tx} dx dt=-2 \int_0^T \int_{L_0}^{L}  \tilde v_{t} \tilde v_{tx} q_x q_t^2 dx dt \nonumber\\&\; \qquad \qquad \qquad\qquad \qquad = \int_0^T \int_{L_0}^{L}  \tilde v_{t}^2 (q_{xx} q_t^2+2q_x q_t q_{tx})dx dt,
		\end{align}
		and
		\begin{align}
			\label{b15}
			\int_0^T \int_{L_0}^{L}  \tilde v_{t}B_3 B_5 \tilde v_{tt} dx dt=&-2 \int_0^T \int_{L_0}^{L}  \tilde v_{t} q_x^2 q_t \tilde v_{tt} dx dt \nonumber\\=& \int_0^T \int_{L_0}^{L}  \tilde v_{t}^2 (q_x^2 q_{tt}+2q_x q_t q_{tx}) dx dt \\&-\int_{L_0}^{L}  \tilde v_{t}^2(T) q_x^2(T) q_t(T) dx+\int_{L_0}^{L}  \tilde v_{t}^2(0) q_x^2(0) q_t(0) dx.\nonumber
		\end{align}
		Moreover, we note that
		\begin{align}
			\label{b16}
			\int_0^T \int_{L_0}^{L}  \tilde v_{t}B_5 B_6 \tilde v dx dt=&2\int_0^T \int_{L_0}^{L}  \tilde v_{t}\tilde v(\rho q_x^4 q_t-q_x^2 q_t^3) dx dt
			\nonumber\\=&-\int_0^T \int_{L_0}^{L} \tilde v^2(\rho q_x^4 q_{tt}+4\rho q_x^3 q_tq_{tx}-2q_x q_t^3q_{tx}-3q_x^2 q_t^2 q_{tt}) dx dt\nonumber\\&+\int_{L_0}^{L}\tilde v^2(T)(\rho q_x^4 (T)q_t(T)-q_x^2(T) q_t^3(T)) dx \nonumber\\&-
			\int_{L_0}^{L}\tilde v^2(0)(\rho q_x^4 (0)q_t(0)-q_x^2(0) q_t^3(0)) dx.
		\end{align}
		Now we take into account that
		\begin{align}
			&q_x=2\tau\mu (x-\tilde L) \eta, \label{q11}\\
			&q_{xx}=2\tau\mu\eta+4\tau \mu^2 (x-\tilde L)^2\eta,  \label{q21}\\
			&q_{xxx}= 12\tau\mu^2 (x-\tilde L)\eta+8\tau \mu^3(x-\tilde L)^3\eta,  \label{q3}\\
			&q_{xxxx}= 12\tau\mu^2\eta+48\tau\mu^3(x-\tilde L)^2\eta+16\tau\mu^4(x-\tilde L)^4\eta, \label{q4}\\
			&q_t=-2\tau\mu(t-\frac{T}{2}) m\eta,  \label{q5}\\
			&q_{tt}=4\tau\mu^2(t-\frac{T}{2})^2 m^2\eta-2\tau\mu m\eta, \label{q6}\\
			&q_{tx}=-4\tau\mu^2 m (x-\tilde L)(t-\frac{T}{2})\eta, \label{q7}\\
			&q_{txx}=-4\tau\mu^2 m (t-\frac{T}{2})\eta-8\tau\mu^3 m (x-\tilde L)^2 (t-\frac{T}{2})  \eta, \label{q8}\\
			&q_{ttx}=8\tau\mu^3 (t-\frac{T}{2})^2(x-\tilde L) m^2\eta-4\tau\mu^2 (x-\tilde L) m\eta, \label{q9}\\
			&q_{ttt}=8\tau \mu^2 m^2(t-\frac{T}{2})\eta -8\tau \mu^3 m^3(t-\frac{T}{2})\eta +4\tau \mu^2 m^2(t-\frac{T}{2})\eta, \\
			&q_{ttxx}=-4\tau \mu^2 m\eta +8\tau \mu^3 m^2 (t-\frac{T}{2})^2 \eta -8\tau \mu^3 m (x-\tilde L)^2 \eta \\&\qquad\qquad\qquad\qquad\qquad\qquad\qquad\qquad+2\tau \mu^4 m^2(x-\tilde L)^2)(t-\frac{T}{2})^2\eta,\\
			&q_{tttx}=16 \tau \mu^3(t\!-\!\frac{T}{2})^2 m^2 (x\!-\!\tilde L)\eta \!-\!16 \tau \mu^4 m^3(t\!-\!\frac{T}{2})^3\eta \!+\!8\tau \mu^3 m^2 (t\!-\!\frac{T}{2})(x\!-\!\tilde L)\eta, \\
			&q_{xxxt}=-24\tau \mu^3 m (x-\tilde L)(t-\frac{T}{2})\eta -16 \tau \mu^4 m (x-\tilde L)^3(t-\frac{T}{2})\eta.
		\end{align}
		and estimate  the last two terms in \eqref{b3} as follows
		\begin{align}
			\label{l1}
			-\rho  \int_0^T \int_{L_0}^{L} \tilde v_{xx}\tilde v_{tt} q_{xx} dx dt&=-\rho \int_0^T \int_{L_0}^{L} \tilde v_{xx}\tilde v_{tt} (2\tau \mu \eta+4\tau \mu^2 (x-\tilde L)^2 \eta) dx dt\nonumber\\&\ge - \rho^2  \int_0^T \int_{L_0}^{L}(2\tau\mu \eta+2\tau\mu^2 (x-\tilde L)^2\eta) \tilde v_{xx}^2dx dt \nonumber\\&\quad- \int_0^T \int_{L_0}^{L} \tilde v_{tt}^2 (\frac{1}{2}\tau \mu\eta+2\tau\mu^2(x-\tilde L)\eta) dx dt
		\end{align}
		and
		\begin{align}
			\label{l2}
			&-4\rho  \int_0^T \int_{L_0}^{L} \tilde v_{xx}\tilde v_{tx} q_{xt} dx dt \nonumber\\&=16\rho \int_0^T \int_{L_0}^{L} \tilde v_{xx}\tilde v_{tx} (\tau m \mu^2 (x-\tilde L)(t-\frac{T}{2}) \eta dx dt\nonumber\\&\ge - \rho^2  \int_0^T\!\! \int_{L_0}^{L} 4\tau\mu^2 (x\!-\!\tilde L)^2 \tilde v_{xx}^2\eta dx dt\!-\! 16\int_0^T \!\!\int_{L_0}^{L} \tilde v_{tt}^2 \tau\mu^2(t\!-\!\frac{T}{2})^2 m^2\eta dx dt.
		\end{align}
		Collecting the second terms in the right-hand sides of \eqref{b8},  \eqref{b11},  and  \eqref{b13} and estimating them from below we arrive at
		\begin{align}
			\label{l3}
			&12\rho  \int_0^T \int_{L_0}^{L} \tilde v_{x}\tilde v_{t} q_{x} q_{t} q_{xx} dx dt \nonumber\\&= -\rho\int_0^T \!\!\int_{L_0}^{L} \tilde v_{x}\tilde v_{t} (48 \tau^3\mu^3 m \mu^2 (x\!-\!\tilde L)(t\!-\!\frac{T}{2}) \eta^3 \!+\!96\tau^3 \mu^4(x\!-\!\tilde L)^3 (t\!-\!\frac{T}{2}) m\eta ^3)dx dt \nonumber\\& \ge -\int_0^T \int_{L_0}^{L} \tilde v_{x}^2( 24\rho^2\tau^3\mu^3(x-\tilde L)^2\eta^3+48\rho^2\tau^3 \mu^4 (x-\tilde L)^4\eta^3) dx dt \nonumber\\&\quad-
			\int_0^T\!\! \int_{L_0}^{L} \tilde v_{t}^2(24\tau^3\mu^3m^2(t\!-\!\frac{T}{2})^2\eta^3\!+\!48\tau^3 \mu^4 m^2 (x\!-\!\tilde L)^2(t\!-\!\frac{T}{2})^2\eta^3) dx dt.
		\end{align}
		Next, using \eqref{q1}--\eqref{q9} we infer that the sum of the term in the right-hand side of \eqref{b1} and the first term in the right-hand side of \eqref{b2}  equals to
		\begin{equation}
			\label{l5}
			\int_0^T \!\!\int_{L_0}^{L} \tilde v_{xx}^2 (3\rho^2\tau\mu\eta+6\rho^2\tau\mu^2(x\!-\!\tilde L)^2 \eta+4\rho\tau\mu^2(t\!-\!\frac{T}{2})^2 m^2 \eta\!-\!2\rho\tau\mu m\eta)dx dt.
		\end{equation}
		Here the relation in brackets can be made positive by choosing $\mu$ large enough. The term in the right-hand side of \eqref{b7} equals to
		\begin{equation}
			\label{l6}
			\int_0^T \int_{L_0}^{L} \tilde v_{tt}^2 (\tau\mu\eta+2\tau\mu^2 (x-\tilde L)^2 \eta)dx dt.
		\end{equation}
		Collecting the first terms in the right-hand side of \eqref{b3}, \eqref{b5}, and \eqref{b6} we obtain
		\begin{equation}
			\label{l7}
			\int_0^T \int_{L_0}^{L} \tilde v_{tx}^2 (2\rho\tau\mu\eta+4\rho\tau\mu^2(x-\tilde L)^2 \eta+4\rho\tau\mu^2(t-\frac{T}{2})^2 m^2 \eta-2\rho\tau\mu m\eta)dx dt,
		\end{equation}
		where the relation in brackets can also be made positive by choosing $\mu$ large enough.
		Analogously, the sum of the first terms in the right-hand sides in \eqref{b4}, \eqref{b9}, \eqref{b10}, \eqref{b13}  equals to
		\begin{multline}
			\label{l8}
			\int_0^T \int_{L_0}^{L} \tilde v_{x}^2 (72\rho^2\tau^3\mu^3(x-\tilde L)^2\eta^3+144\rho^2\tau^3\mu^4(x-\tilde L)^4 \eta^3\\+48\tau^3\mu^3(t-\frac{T}{2})^4 m^4 \eta^3-24\tau^3\mu^3 m^3(t-\frac{T}{2})^2\eta^3)dx dt
		\end{multline}
		and the sum of the first terms in \eqref{b8}, \eqref{b11}, \eqref{b14}, \eqref{b15} equals to
		\begin{multline}
			\label{l9}
			\int_0^T \int_{L_0}^{L} \tilde v_{t}^2 (48\rho\tau^3\mu^3(x-\tilde L)^2\eta^3+96\rho\tau^3\mu^4(x-\tilde L)^4 \eta^3\\+48\tau^3\mu^4(t-\frac{T}{2})^2(x-\tilde L)^2 m^2 \eta^3-8\tau^3\mu^3 m (x-\tilde L)^2\eta^3)dx dt,
		\end{multline}
		where the relations in brackets can also be made positive by choosing $\mu$ large enough. \par
		We note that
		\begin{align}
			\label{l10}
			&\int_0^T \int_{L_0}^{L} \tilde v^2(\frac{15}{2} \rho^2 q_x^4 q_{xx}- \rho q_x^4 q_{tt}-4\rho q_x^3 q_tq_{tx})dx dt \nonumber\\&\ge 160\rho^2\int_0^T \int_{L_0}^{L} \tilde v^2  \tau^5\mu^6(x-\tilde L)^6\eta^5 dx dt.
		\end{align}
		This term can be made large enough by choosing $\tau$ so that the sum of all terms containing $\tilde v^2$ is positive.
		Collecting \eqref{i2}--\eqref{l10} with $\mu$ and $\tau$ large enough we obtain that there exist $C_1, C_2>0$ such that
		\begin{align}
			\label{car1}
			&\int_0^T \int_{L_0}^{L} I_1 I_2 dx dt \nonumber\\&\ge C_1\bigg[\tau\mu^2\int_0^T \int_{L_0}^{L}\eta (\tilde v_{xx}^2+\tilde v_{tt}^2+v_{xt}^2) dx dt+\tau^3\mu^4\int_0^T \int_{L_0}^{L}\eta^3 (\tilde v_{x}^2+\tilde v_{t}^2) dx dt\nonumber\\&\quad+\tau^5\mu^6\int_0^T\!\! \int_{L_0}^{L}\eta^5 \tilde v^2 dx dt\bigg]\!-\!C_2[\mu\tau\eta(T)E_{v1}(T)\!+\!\mu\tau\eta(0)E_{v1}(0)\!+\!\mu^3\tau^3\eta^3(T)E_{v2}(T)\nonumber\\&\quad+\mu^3\tau^3\eta^3(0)E_{v2}(0)+\mu^5\tau^5\eta^5(T)E_{v3}(T)+\mu^5\tau^5\eta^5(0)E_{v3}(0)],
		\end{align}
		where
		\begin{gather}
			E_{v1}(t)= \int_{L_0}^{L} v_{xx}^2(t) dx+\int_{L_0}^{L} v_{tx}^2(t) dx\label{car21}\\
			E_{v2}(t)= \int_{L_0}^{L} v_{x}^2(t) dx+\int_{L_0}^{L} v_{t}^2(t) dx\label{car22}\\
			E_{v3}(t)= \int_{L_0}^{L} v^2(t) dx.\label{car23}
		\end{gather}
		It follows from \eqref{per}, \eqref{1i}-- \eqref{1s}, \eqref{q1}-- \eqref{q9} that
		\begin{align}
			\label{car3}
			&2 \int_0^T \int_{L_0}^{L} I_1I_2 dx dt \nonumber\\&\le \int_0^T \int_{L_0}^{L} |I_1+I_2|^2 dx dt\le  C\left(\int_0^T \int_{L_0}^{L}\theta^2 |P v|^2 dx dt
			+\int_0^T \int_{L_0}^{L} |S v|^2 dx dt\right)\\&\le
			C\bigg(\int_0^T \!\!\int_{L_0}^{L}\theta^2 |P v|^2 dx dt
			\!+\!\tau^2\mu^4\int_0^T \int_{L_0}^{L}( \tilde v_x^2\!+\!\tilde v_t^2) \eta^2 dx dt \!+\!\tau^4\mu^6\int_0^T\!\! \int_{L_0}^{L} \tilde v^2 \eta^2 dx dt\bigg).\nonumber
		\end{align}
		Collecting \eqref{car1}--\eqref{car3}, and choosing $\tau$ large enough we obtain the following Carleman estimate
		\begin{align}
			\label{car4}
			&\int_0^T \int_{L_0}^{L}\theta^2 |P v|^2 dx dt \nonumber\\&\ge  C_1\bigg[\tau\mu^2\int_0^T \int_{L_0}^{L}\eta \theta^2 ( v_{xx}^2+\tilde v_{tt}^2+v_{xt}^2) dx dt \nonumber\\&\quad+\tau^3\mu^4\int_0^T \int_{L_0}^{L}\eta^3\theta^2 ( v_{x}^2+ v_{t}^2) dx dt+\tau^5\mu^6\int_0^T \int_{L_0}^{L}\eta^5 \theta^2\tilde v^2 dx dt\bigg]\nonumber\\&\quad-C_2\bigg[\mu\tau\eta(T)E_{v1}(T)+\mu\tau\eta(0)E_{v1}(0)+\mu^3\tau^3\eta^3(T)E_{v2}(T)+\mu^3\tau^3\eta^3(0)E_{v2}(0)\nonumber\\&\quad+\mu^5\tau^5\eta^5(T)E_{v3}(T)+\mu^5\tau^5\eta^5(0)E_{v3}(0) \bigg].
		\end{align}
		The Carleman estimate (see the ideas of the proof in e.g.    \cite{Bod} and references therein) for the operator $Q=\partial_t^2-\rho \partial_x^2$ for a function
		\begin{equation*}
			\hat u\in H^2((L_0, L)\times (0, T)), \;\hat u(L_0)=\hat u_x(L_0)=\hat u(L)=0
		\end{equation*}
		is as follows
		\begin{align}
			\label{car5}
			&\int_0^T \int_{L_0}^{L}\theta^2 |Q \hat u|^2 dx dt \nonumber\\&\ge  C_3\left[\tau\mu\int_0^T \int_{L_0}^{L}\eta \theta^2( \hat u_{x}^2+\hat u_{t}^2) dx dt+\tau^3\mu^4\int_0^T \int_{L_0}^{L}\eta^3\theta^2 \hat u^2 dx dt\right] \\&\quad-C_4\left[\mu\tau \eta(T) E_{\hat u2}(T)+\mu\tau\eta(0) E_{\hat u2}(0)+\mu^3\tau^3\eta^3(T) E_{\hat u3}(T)+\mu^3\tau^3\eta^3(0) E_{\hat u3}(0) \right],\nonumber
		\end{align}
		Consequently, for a function
		\begin{equation*}
			u\in H^4((L_0, L)\times (0, T)), \; u_{xx}(L_0)=u_{xxx}(L_0)=u_{xx}(L)=0
		\end{equation*}
		we have from \eqref{car5} choosing $\hat u=u_{xx}$
		\begin{align}
			\label{car6}
			&\int_0^T \int_{L_0}^{L}\theta^2 |Q u_{xx}|^2 dx dt \nonumber\\&\ge  C_3\left[\tau\mu\int_0^T \int_{L_0}^{L}\eta \theta^2( u_{xxx}^2+ u_{txx}^2) dx dt+\tau^3\mu^4\int_0^T \int_{L_0}^{L}\eta^3\theta^2 u_{xx}^2 dx dt\right]\\&\quad-C_4\left[\mu\tau\eta(T)E_{u4}(T)\!+\!\mu\tau\eta(0)E_{u4}(0)+\mu^3\tau^3\eta^3(T)E_{u1}(T)+\mu^3\tau^3\eta^3(0)E_{u1}(0) \right],\nonumber
		\end{align}
		where
		\begin{equation}
			E_{v4}(t)= \int_{L_0}^{L} v_{xxx}^2(t) dx+\int_{L_0}^{L} v_{txx}^2(t) dx\label{car41}.
		\end{equation}
		\item[{\it Step 3. Observability inequality.}]
		From estimates \eqref{car5}, \eqref{car6} and equations \eqref{du} and
		\begin{equation}
			\label{dvv}
			\beta_2 v_{dttx}- \left(v_{dx}+1/2u_{dx}(u_x(t+h)+u_x(t))\right)_{xx}=0
		\end{equation}
		we infer the following estimate
		\begin{align}
			\label{ccc1}
			&\tau\mu\int_0^T \int_{L_0}^{L}\eta \theta^2( u_{dxxx}^2+ u_{dtxx}^2) dx dt+\tau^3\mu^4\int_0^T \int_{L_0}^{L}\eta^3\theta^2 u_{dxx}^2 dx \nonumber\\&\quad +\tau\mu^2\int_0^T \int_{L_0}^{L}\eta \theta^2 ( v_{dxx}^2+\tilde v_{dtt}^2+v_{dxt}^2) dx dt+\tau^3\mu^4\int_0^T \int_{L_0}^{L}\eta^3\theta^2 ( v_{dx}^2+ v_{dt}^2) dx dt \nonumber\\&\quad+\tau^5\mu^6\int_0^T \int_{L_0}^{L}\eta^5 \theta^2\tilde v_d^2 dx dt-C\bigg[\mu\tau\eta(T)E_{v_d1}(T)+\mu\tau\eta(0)E_{v_d1}(0) \\&\quad+\mu^3\tau^3\eta^3(T)E_{v_d2}(T)\!+\!\mu^3\tau^3\eta^3(0)E_{v_d2}(0)\!+\!\mu^5\tau^5\eta^5(T)E_{v_d3}(T)\!+\!\mu^5\tau^5\eta^5(0)E_{v_d3}(0) \nonumber\\&\quad+\mu^3\tau^3\eta^3(T)E_{u_d1}(T)\!+\!\mu^3\tau^3\eta^3(0)E_{u_d1}(0)\!+\!\mu\tau\eta(T)E_{u_d4}(T)\!+\!\mu\tau\eta(0)E_{u_d4}(0)\bigg]\!\le\! 0 \nonumber
		\end{align}
		choosing in \eqref{ccc1} $\tau$ and $\mu$ large enough and taking into account that $E_{v_d3}(t)\le C E_{v_d2}(t)$ we arrive at
		\begin{equation}
			\label{cc1}
			\tau\mu\int_{t_0}^{t_1} \int_{L_0}^{L}E_{u_d, v_d}(t) dx dt -C\mu^5\tau^5e^{-\sigma_0 \mu}\left[E_{u_d, v_d}(0)+E_{u_d, v_d}(T)\right]\le 0,
		\end{equation}
		where
		\[E_{u_d,v_d}(t)=E_{v_d1}(t)+E_{v_d2}(t)+E_{u_d1}(t)+E_{u_d2}(t)+E_{u_d4}(t).\]
		We multiply equation \eqref{du} by $- u_{dtxx}$ and equation \eqref{dv} by $- v_{dtxx}$ and integrate by parts over the intervals $[L_0, L]$ and
		$[s,t]\subset[0,T]$ to get
		\begin{align}
			\label{e4e}
			&E_{u_d4}(t)+E_{u_d1}(t)+E_{v_d1}(t) \nonumber\\&\le E_{u_d4}(s)\!+\!E_{u_d1}(s)\!+\!E_{v_d1}(s) \!+\!C_T\int_{s}^{t}(E_{u_d4}(\xi)\!+\!E_{u_d1}(\xi)\!+\!E_{v_d1}(\xi))d\xi.
		\end{align}
		Analogously, after multiplication of equation \eqref{du} by $u_{dt}$ and equation \eqref{dv} by $v_{dt}$ and integration by parts over the intervals $[L_0, L]$ and
		$[s,t]\subset[0,T]$  we arrive at
		\begin{align}
			\label{e2e}
			&E_{u_d1}(t)+E_{u_d2}(t)+E_{v_d2}(t) \nonumber\\&\le E_{u_d1}(s)\!+\!E_{u_d2}(s)\!+\!E_{v_d2}(s) \!+\!C_T\int_{s}^{t}(E_{u_d1}(\xi)\!+\!E_{u_d2}(\xi)\!+\!E_{v_d2}(\xi))d\xi.
		\end{align}
		It follows from \eqref{e4e}, \eqref{e2e}, and  the Gronwall's lemma that for $0\le s\le t\le T$
		\begin{equation}
			E_{u_d, v_d}(t)\le E_{u_d, v_d}(s)e^{C_T(t-s)}\label{es1}.
		\end{equation}
		Analogously,
		\begin{equation}
			E_{u_d, v_d}(s)\le E_{u_d, v_d}(t)e^{C_T(t-s)}.\label{es2}
		\end{equation}
		Choosing $t=T$, $s=t$ in \eqref{es1} and $s=0$ in \eqref{es2} and summing up the results we obtain
		\begin{equation}
			\label{ess}
			E_{u_d, v_d}(T)+E_{u_d, v_d}(0)\le E_{u_d, v_d}(t)e^{C_T T}.
		\end{equation}
		Substituting \eqref{ess} into \eqref{cc1} and taking into account \eqref{q2} we arrive at
		\begin{equation}
			\label{finobs}
			(\tau\mu(t_1-t_0)e^{\sigma_2\mu-C_T T} -C\mu^5\tau^5e^{-\sigma_1 \mu})\left[E_{u_d, v_d}(0)+E_{u_d, v_d}(T)\right]\le 0,
		\end{equation}
		where $\mu$ can be chosen large enough, so that the constant $\tau\mu(t_1-t_0)e^{\sigma_2\mu-C_T T} -C\mu^5\tau^5e^{-\sigma_1 \mu}$
		is positive.This immediately gives $u_d=v_d=0$ for any $t\ge 0$, consequently, the system $(S_t,H)$ is gradient.
	\end{proof}
	Now we state our main result.
	\begin{theorem}
		Let assumptions of Theorem 3, Theorem 4, and Theorem 5 hold true. Moreover, let
		\begin{equation}
			g_2(x)=g_4(x)=0.
		\end{equation}
		Then, the dynamical system $(S_t, H)$ generated by \eqref{1}-\eqref{24} possesses a compact global attractor possessing properties \eqref{conv-N}, \eqref{7.4.1}.
	\end{theorem}
	\begin{proof}
		In view  of Theorem 2, Theorem 4, and Theorem 5 our remaining task is to show the boundedness of the set of stationary points and the set $W_R=\{Z: \EuScript L(Z)\le R\}$, where $\EuScript L$ is given by \eqref{lap}.\par
		The second statement follows immediately from the structure of function $\EuScript L$ and Lemma 1.\par
		The first statement can be easily shown by the substitution of $ \Psi=(\phi, u, \omega, v)$ into \label{sol_def} and application of energy-like estimates and Lemma 1 for stationary solutions.
	\end{proof}
	\begin{remark}
		From the point of view of applications  it is interesting to consider the system
		\begin{align}
			&\beta_1 \phi_{tt}-\mu_1 \phi_{ttxx}-\kappa\phi_{txx} +\lambda_1
			\phi_{xxxx}-\delta_1\left(\left[\phi_x\left(\omega_x+1/2\phi_x^2\right) \right]_x \right)=g_1(x,t),\\
			&\rho_1 \omega_{tt}+ \gamma \omega_{t}- \delta_1\left(\omega_x+1/2\phi_x^2\right)_x=g_2(x,t).\qquad\qquad\;\;\;t>0,\;\;x\in (0, L_0)
			\\
			&\beta_2 u_{tt}-\mu_2 u_{ttxx}+\lambda_2
			u_{xxxx}-\delta_2\left(\left[u_x\left(v_x+1/2u_x^2\right) \right]_x \right)=g_3(x,t),\\
			&\rho_2 v_{tt}- \delta_2\left(v_x+1/2u_x^2\right)_x=g_4(x,t),\qquad\qquad\qquad\;\;\;\;\;\;t>0,\;\;x\in ( L_0, L)
		\end{align}
		with the transmission boundary condition
		\begin{equation}
			\delta_1(\omega_x+1/2\phi_x^2)(L_0, t)=\delta_2 (v_x+1/2u_x^2)(L_0, t)
		\end{equation}
		instead of \eqref{7}. This means that all the physical properties of two parts of the beam are different. However, in this case one can prove only the existence of weak solutions and this is an open question how to show the higher order estimates without use of strong solutions in the arguments.
	\end{remark}
	\section*{Acknowledgements}
	The author would like to thank the Referee whose valuable
	comments and suggestions helped to improve the present work.\\
	The author was partially supported by successively the Volkswagen Foundation grant within the
	frameworks of the international project “ From Modeling and Analysis to Approximation” and the Volkswagen Foundation grant for the project  ``Dynamic problems in elasticity" at Humboldt-Universit\"at zu Berlin, Funding for Refugee Scholars and Scientists from Ukraine.

	\medskip

\end{document}